\documentclass[article,reqno]{amsart}

\usepackage{amssymb}

\newtheorem{theorem}{Theorem}
\theoremstyle{plain}

\newtheorem{corollary}{Corollary}

\newtheorem{lemma}{Lemma}

\newtheorem{proposition}{Proposition}
\newtheorem{remark}{Remark}

\numberwithin{equation}{section}

\begin{document}
\title[SILT of the $(\alpha ,d,\beta )$-superprocess]{Self-Intersection
Local Time of $(\alpha ,d,\beta )$-superprocess}
\author{L. Mytnik$\hbox{}^{1}$}
\address{Technion-Israel Institute of Technology\\
Faculty of Industrial Engineering and Management\\
Haifa 32000, Israel}
\email{leonid@ie.technion.ac.il}
\author{J. Villa $\hbox{}^{2}$}
\address{Universidad Aut\'{o}noma de Aguascalientes\\
Departamento de Matem\'{a}ticas y F\'{\i}sica\\
Av. Universidad 940 C.P. 20100, Aguascalientes, Ags., M\'{e}xico}
\email{jvilla@correo.uaa.mx}
\thanks{(2) Thanks go to UAA for allow this postdoctoral stay.}
\thanks{(1), (2) Research supported in part by the Israel Science Foundation.%
}

\begin{abstract}
The existence of self-intersection local time (SILT), when the time diagonal
is intersected, of the $(\alpha ,d,\beta )$-superprocess is proved for $%
d/2<\alpha $\ and for a renormalized SILT when $d/(2+(1+\beta )^{-1})<\alpha
\leq d/2$. We also establish Tanaka-like formula for SILT.
\end{abstract}

\maketitle

\begin{center}

\end{center}

\section{Introduction and statement of results}

This paper is devoted to the proof of existence of self-intersection local
time (SILT) of $(\alpha ,d,\beta )$-superprocesses for $0<\beta <1$. Let us
introduce some notation. Let $\mathcal{B}_{b}(\mathbb{R}^{d})$ (respectively 
$\mathcal{C}_{b}(\mathbb{R}^{d})$) be the family of all bounded
(respectively, bounded continuous) Borel measurable functions on $\mathbb{R}%
^{d}$, and $\mathcal{M}_{F}(\mathbb{R}^{d})$ be the set of all finite Borel
measures on $\mathbb{R}^{d}$. The integral of a function $f$ with respect to
a measure $\mu $ is denoted by $\mu (f)$. If $E$ is a metric space we denote
by $D([0,+\infty ),E)$ the space of all c\`{a}dl\`{a}g $E$-valued paths with
the Skorohod topology. We will use $c$ to denote a positive and finite
constant whose value may vary from place to place. A constant of the form $%
c(a,b,...)$ means that this constant depends on parameters $a,b,...$.

Let $(\Omega ^{\prime },\mathcal{F}^{\prime },\mathcal{F}_{\cdot }^{\prime
},P^{\prime })$ be a filtered probability space where the $(\alpha ,d,\beta
) $-superprocess $X=\{X_{t}:t\geq 0\}$ is defined. That is, by $X$ we mean a 
$\mathcal{M}_{F}(\mathbb{R}^{d})$-valued, time homogeneous, strong Markov
process with c\`{a}dl\`{a}g sample paths, such that for any non-negative
function $\varphi \in \mathcal{B}_{b}(\mathbb{R}^{d})$, 
\begin{equation*}
E\left[ \left. \exp \left( -X_{t}(\varphi )\right) \right| X_{0}=\mu \right]
=\exp \left( -\mu (V_{t}(\varphi ))\right) ,
\end{equation*}
where $\mu \in \mathcal{M}_{F}(\mathbb{R}^{d})$ and $V_{t}(\varphi )$
denotes the unique non-negative solution of the following evolution equation 
\begin{equation*}
v_{t}=S_{t}\varphi -\int_{0}^{t}S_{t-s}\left( (v_{s})^{1+\beta }\right) ds,%
\text{ \ \ }t\geq 0.
\end{equation*}
Here $\{S_{t}:t\geq 0\}$ denotes the semigroup corresponding to the
fractional Laplacian operator $\Delta _{\alpha }$.

Another way to characterize the $(\alpha ,d,\beta )$-superprocess $X$ is by
means of the following martingale problem: 
\begin{equation}
\left\{ 
\begin{tabular}{l}
For all $\varphi \in D(\Delta _{\alpha })$ (domain of $\Delta _{\alpha }$)
and $\mu \in \mathcal{M}_{F}(\mathbb{R}^{d})$, \\ 
$X_{0}=\mu $, and $M_{t}(\varphi )=X_{t}(\varphi )-X_{0}(\varphi
)-\int_{0}^{t}X_{s}(\Delta _{\alpha }\varphi )ds$, \\ 
is a $\mathcal{F}_{t}^{\prime }$-martingale.%
\end{tabular}%
\ \right.   \label{mp}
\end{equation}%
If $\beta =1$ then $M_{\cdot }(\varphi )$ is a continuous martingale. In
this paper we are interested in the case of $0<\beta <1,$ and here $%
M_{t}(\varphi )$ is a purely discontinuous martingale. This martingale can
be expressed as 
\begin{equation}
M_{t}(\varphi )=\int_{0}^{t}\int_{\mathbb{R}^{d}}\varphi (x)M(ds,dx),
\label{storepmtg}
\end{equation}%
where $M(ds,dx)$ is a martingale measure, it and the stochastic integral
with respect to such martingale measure is defined in \cite{L-M} (or in
Section II.3 of \cite{I-W}).

The SILT is heuristically defined by 
\begin{equation*}
\gamma _{X}(B)=\int_{B}\int_{\mathbb{R}^{2d}}\delta
(x-y)X_{s}(dx)X_{t}(dy)dsdt,
\end{equation*}%
where $B\subset \lbrack 0,\infty )\times \lbrack 0,\infty )$ is a bounded
Borel set and $\delta $ is the Dirac delta function. Let $D=\{(t,t):t>0\}$
be the time diagonal on $\mathbb{R}_{+}\times \mathbb{R}_{+}$. For $\beta =1$%
, $B\cap D=\varnothing $ and $d\leq 7$, Dynkin \cite{D} proved the existence
of SILT, $\gamma _{X}$, for a very general class of continuous
superprocesses. Also, from the Dynkin's works follows the existence of SILT
when $\beta =1$, $B\cap D\neq \varnothing $ and $d\leq 3$ (see \cite{A}).
For $\beta =1$, $d=4,5$ and $B\cap D\neq \varnothing ,$ Rosen \cite{R}
proved the existence of a renormalized SILT for the $(\alpha ,d,1)$%
-superprocess. A Tanaka-like formula for the local time of $(\alpha ,d,2)$%
-superprocess was established by Adler and Lewin in \cite{A-L2}. The same
authors derived a Tanaka-like formula for self-intersection local time for $%
(\alpha ,d,2)$-superprocess (see \cite{A-L1}). In this paper we are going to
extend the above results for the case of $0<\beta <1$.

The usual way to give a rigorous definition of SILT is to take a sequence $%
(\varphi _{\varepsilon })_{\varepsilon >0}$ of smooth functions that
converges in distribution to $\delta $, define the approximating SILTs 
\begin{equation*}
\gamma _{X,\varepsilon }(B)=\int_{B}\int_{\mathbb{R}^{2d}}\varphi
_{\varepsilon }(x-y)X_{s}(dx)X_{t}(dy)dsdt,\text{ \ \ }B\subset \mathbb{R}%
_{+}\times \mathbb{R}_{+},
\end{equation*}%
and prove that $(\gamma _{X,\varepsilon }(B))_{\varepsilon >0}$ converges,
in some sense (it is usually taken $L^{2}(P^{\prime }),$ $L^{1+\beta
}(P^{\prime })$, $L^{1}(P^{\prime })$, distribution or in probability), to a
random variable $\gamma _{X}(B)$. In what follows we choose $\varphi
_{\varepsilon }=p_{\varepsilon }$, where $p_{\varepsilon }$ is the $\alpha $%
-stable density, given by 
\begin{equation*}
p_{\varepsilon }(x,y)=\frac{1}{(2\pi )^{d}}\int_{\mathbb{R}^{d}}e^{-i(z\cdot
(x-y))-\varepsilon |z|^{\alpha }}dz,\text{ \ \ }x,y\in \mathbb{R}^{d},
\end{equation*}%
when $0<\alpha <2$ and 
\begin{equation*}
p_{\varepsilon }(x,y)=\frac{1}{(2\pi \varepsilon )^{d/2}}e^{-|x-y|^{2}/2%
\varepsilon },\text{ \ \ }x,y\in \mathbb{R}^{d},
\end{equation*}%
for $\alpha =2$. In this paper we will consider the particular case when $%
B=\{(t,s):0\leq s\leq t\leq T\}$. Here we denote $\gamma _{X,\varepsilon
}(B) $ by $\gamma _{X,\varepsilon }(T)$, that is 
\begin{equation*}
\gamma _{X,\varepsilon }(T)=\int_{0}^{T}\int_{0}^{t}\int_{\mathbb{R}%
^{2d}}p_{\varepsilon }(x-y)X_{s}(dx)X_{t}(dy)dsdt,\;\;\forall T\geq 0.
\end{equation*}%
Moreover, we are going to consider the renormalized SILT 
\begin{equation*}
\tilde{\gamma}_{X,\varepsilon }(T)=\gamma _{X,\varepsilon }(T)-e^{\lambda
\varepsilon }\int_{0}^{T}\int_{\mathbb{R}^{d}}X_{s}(G^{\lambda ,\varepsilon
}(x-\cdot ))X_{s}(dx)ds,
\end{equation*}%
where 
\begin{equation*}
G^{\lambda ,\varepsilon }(x)=\int_{\varepsilon }^{\infty }e^{-\lambda
t}p_{t}(x)dt,\text{ \ }\lambda ,\varepsilon \geq 0.
\end{equation*}

Notice that $G^{\lambda ,\varepsilon }(x)\uparrow G^{\lambda
,0}(x)=G^{\lambda }(x)$ as $\varepsilon \downarrow 0$ for all $x\in \mathbb{R%
}^{d}$, and if $d>\alpha $ then 
\begin{equation}
G(x)=G^{0,0}(x)=c(\alpha ,d)\left\vert x\right\vert ^{\alpha -d},
\label{expfungreen}
\end{equation}%
where 
\begin{equation*}
c(\alpha ,d)=\frac{\Gamma ((d-\alpha )/2)}{2^{\alpha /2}\pi ^{d/2}\Gamma
(\alpha /2)},
\end{equation*}%
and $\Gamma $ is the usual Gamma function. $G$ is called Green function of $%
\Delta _{\alpha }$. Also notice that, for $\lambda >0$, we have 
\begin{equation}
\Delta _{\alpha }G^{\lambda ,\varepsilon }(x)=\lambda G^{\lambda
,\varepsilon }(x)-e^{-\lambda \varepsilon }p_{\varepsilon }(x),\text{ \ }%
x\in \mathbb{R}^{d}.  \label{genagreeb}
\end{equation}

\bigskip

\noindent Now we are ready to present our main result.

\begin{theorem}
\label{Theomain}Let $X$ be the $(\alpha ,d,\beta )$-superprocess with
initial measure $X_{0}(dx)=\mu (dx)=h(x)dx$, where $h$ is bounded and
integrable with respect to Lebesgue measure on $\mathbb{R}^{d}$. Let $M$ be
the martingale measure which appears in the martingale problem~(\ref{mp})
for $X$.

\begin{itemize}
\item[\textbf{(a)}] Let $d/2<\alpha $. Then there exists a process $\gamma
_{X}=\{\gamma _{X}(T):T\geq 0\}$ such that for every $T>0$, $\delta>0$ 
\begin{equation*}
P\left( \sup_{t\leq T}\left| \gamma _{X,\varepsilon }(t)-\gamma
_{X}(t)\right|>\delta \right) \rightarrow 0,\;as\;\varepsilon \downarrow 0.
\end{equation*}
Moreover, for any $\lambda >0,$ 
\begin{eqnarray}
\gamma _{X}(T) &=&\lambda \int_{0}^{T}\int_{0}^{t}\int_{\mathbb{R}%
^{2d}}G^{\lambda }(x-y)X_{s}(dx)X_{t}(dy)dtds  \label{tkfmala} \\
&&-\int_{0}^{T}\int_{\mathbb{R}^{d}}X_{T}(G^{\lambda }(x-\cdot ))X_{s}(dx)ds
\notag \\
&&+\int_{0}^{T}\int_{\mathbb{R}^{d}}X_{s}(G^{\lambda }(x-\cdot ))X_{s}(dx)ds
\notag \\
&&+\int_{0}^{T}\int_{\mathbb{R}^{d}}\int_{0}^{t}\int_{\mathbb{R}%
^{d}}G^{\lambda }(x-y)M(ds,dy)X_{t}(dx)dt, \;\;\;\mathrm{a.s.}  \notag
\end{eqnarray}

\item[\textbf{(b)}] Let $d/(2+(1+\beta )^{-1})<\alpha \leq d/2$. Then there
exists a process $\tilde{\gamma}_{X}=\{\tilde{\gamma}_{X}(T):T\geq 0\}$ such
that for every $T>0$, $\delta>0$ 
\begin{equation*}
P\left( \sup_{t\leq T}\left| \tilde{\gamma}_{X,\varepsilon }(t)-\tilde{\gamma%
}_{X}(t)\right| >\delta\right) \rightarrow 0,\;as\;\varepsilon \downarrow 0.
\end{equation*}
Moreover, for any $\lambda >0,$ 
\begin{eqnarray}
\tilde{\gamma}_{X}(T) &=&\lambda \int_{0}^{T}\int_{0}^{t}\int_{\mathbb{R}%
^{2d}}G^{\lambda }(x-y)X_{s}(dx)X_{t}(dy)dtds  \label{tkfmalarenor} \\
&&-\int_{0}^{T}\int_{\mathbb{R}^{d}}X_{T}(G^{\lambda }(x-\cdot ))X_{s}(dx)ds
\notag \\
&&+\int_{0}^{T}\int_{\mathbb{R}^{d}}\int_{0}^{t}\int_{\mathbb{R}%
^{d}}G^{\lambda }(x-y)M(ds,dy)X_{t}(dx)dt,\;\;\mathrm{a.s.}  \notag
\end{eqnarray}
\end{itemize}

The processes $\gamma _{X}$ and $\tilde{\gamma}_{X}$ are called SILT and
renormalized SILT of $X$, respectively, and (\ref{tkfmala}) and (\ref%
{tkfmalarenor}) are called Tanaka-like formula for SILT.
\end{theorem}

\begin{remark}
It is interesting to note that our bound on dimensions 
\begin{equation*}
d<(2+(1+\beta )^{-1})\alpha ,
\end{equation*}%
for renormalized SILT does not converge, as $\beta \uparrow 1$, to the bound 
$d<3\alpha $ established by Rosen~\cite{R} for finite variance superprocess (%
$\beta =1$). For example, simple algebra shows that for $5/3<\alpha <2$,
there is a SILT for finite variance superprocees in dimensions $d\leq 5$.
However for any $\beta \in ((3\alpha -5)/(5-2\alpha ),1)$ we get the
exixtence of SILT only in dimensions $d\leq 4$ and this bound not improve to 
$d\leq 5$ if $\beta \uparrow 1$. So, our bound is more restrictive, and we
believe that it is related to the fact that $(\alpha ,d,\beta )$%
-superprocess (with $\beta <1$) has jumps. Our conjecture is that for $\beta
<1$, the renormalized SILT defined by~(\ref{tkfmalarenor}) does not exist in
dimensions greater than $(2+(1+\beta )^{-1})\alpha $.
\end{remark}

The common ways to prove the existence of SILT for the finite variance
superprocesses (see e.g. \cite{A-L1}, \cite{R}) do not work here. The reason
for this is that such proofs strongly rely on the existence of high moments
of $X$ (at least of order four), and $(\alpha ,d,\beta )$-superprocess $X$
has moments of order less than $1+\beta $. To overcome this difficulty let
us consider the path properties of $X$ more carefully. It is well known (see
Theorem 6.1.3 of \cite{Da}) that, for $0<\beta <1,$ the $(\alpha ,d,\beta )$%
-su\-per\-pro\-cess $X$ is a.s. discontinuous and has jumps of the form $%
\Delta X_{t}=r\delta _{x}$, for some $r>0,$ $x\in {\mathbb{R}}^{d}$. Here $%
\delta _{x}$ denotes the Dirac measure concentrated at $x$. Let 
\begin{equation}
N_{X}(dx,dr,ds)=\sum_{\{(x,r,s):\Delta X_{s}=r\delta _{x}\}}\delta
_{(x,r,s)},  \label{dmc}
\end{equation}%
be a random point measure on $\mathbb{R}^{d}\times \mathbb{R}_{+}\times 
\mathbb{R}_{+}$ with compensator measure $\hat{N}_{X}$ given by 
\begin{equation}
\hat{N}_{X}(dx,dr,ds)=\eta r^{-2-\beta }drX_{s}(dx)ds,  \label{cmc}
\end{equation}%
where 
\begin{equation*}
\eta =\frac{\beta (\beta +1)}{\Gamma (1-\beta )}.
\end{equation*}%
Let $K>0$ fix$.$ From \cite{K-R} and \cite{Da} we have that the $(\alpha
,d,\beta )$-superprocess $X$ has the following decomposition: Let $\varphi
\in D(\Delta _{\alpha })$, $t\geq 0$, 
\begin{eqnarray}
X_{t}(\varphi ) &=&\mu (\varphi )+\int_{0}^{t}X_{s}(\Delta _{\alpha }\varphi
)ds-C_{\beta }(K)\int_{0}^{t}X_{s}(\varphi )ds  \notag \\
&&+\int_{0}^{t}\int_{0}^{K}\int_{\mathbb{R}^{d}}r\varphi (x)\widetilde{N}%
_{X}(dx,dr,ds)  \notag \\
&&+\int_{0}^{t}\int_{K}^{\infty }\int_{\mathbb{R}^{d}}r\varphi
(x)N_{X}(dx,dr,ds),  \label{bigjumps}
\end{eqnarray}%
where $\widetilde{N}_{X}=N_{X}-\hat{N}_{X}$ is a martingale measure and 
\begin{equation*}
C_{\beta }(K)=\frac{\eta }{\beta K^{\beta }}.
\end{equation*}%
As we have mentioned already, one of the problems of working with the $%
(\alpha ,d,\beta )$-superprocess $X$ is dealing with \textquotedblleft
big\textquotedblright\ jumps. In fact, the \textquotedblleft
big\textquotedblright\ jumps produce the infinite variance of the process
and they appear in the term corresponding to the integral with respect to $%
N_{X}$ on\ (\ref{bigjumps}). So, the first step in the establishing the
existence of SILT for $(\alpha ,d,\beta )$-superprocess $X$ is to
\textquotedblleft eliminate\textquotedblright\ those jumps. This is achieved
via introducing the following auxiliary process.

Let us considerer the canonical space, $\Omega ^{\circ }=D([0,\infty ),M_{F}(%
\mathbb{R}^{d}))$, $\mathcal{F}^{\circ }=\mathcal{B}(\Omega ^{\circ })$ and $%
\mathcal{F}_{t}^{\circ }=\sigma \{Y_{r}^{K}:0\leq r\leq t\}$, where $%
Y_{r}^{K}(\omega ^{\circ })=\omega ^{\circ }(r)$. For any $\mu \in M_{F}(%
\mathbb{R}^{d})$ there exists (see \cite{Da}) a measure $Q_{\mu }$ on $%
(\Omega ^{\circ },\mathcal{F}^{\circ })$, such that for any non-negative
function $\varphi \in \mathcal{B}_{b}(\mathbb{R}^{d})$%
\begin{equation}
E_{\mu }\left[ \text{exp}\left( -Y_{t}^{K}(\varphi )\right) |\mathcal{F}%
_{s}^{\circ }\right] =\text{exp}\left( -Y_{s}^{K}\left( V_{t-s}^{K}(\varphi
)\right) \right) ,\;\forall 0<s\leq t,  \label{equt:3}
\end{equation}%
and 
\begin{equation}
E_{\mu }\left[ \text{exp}\left( -Y_{t}^{K}(\varphi )\right) \right] =\text{%
exp}\left( -\mu \left( V_{t}^{K}(\varphi )\right) \right) ,\;\;\;\forall t>0,
\label{equat1pp}
\end{equation}%
(notice that the expectation is taken here with respect to the measure $%
Q_{\mu }$). $V_{t}^{K}$ is the unique non-negative solution for the
non-linear equation%
\begin{equation}
\left\{ 
\begin{tabular}{rll}
$\frac{\partial v_{t}^{K}}{\partial t}$ & $=$ & $(\Delta _{\alpha }-C_{\beta
}(K))v_{t}^{K}-\Phi ^{K}(v_{t}^{K}),$ \\ 
$v_{0}^{K}$ & $=$ & $\varphi ,$%
\end{tabular}%
\right.   \label{equt:2}
\end{equation}%
where 
\begin{equation}
\Phi ^{K}(x)=\eta \int_{0}^{K}\left( e^{-ux}-1+ux\right) u^{-\beta -2}du.
\label{defapxintpsim}
\end{equation}%
Note that when $K=\infty $ the resulting process $Y^{\infty }$ and the
regular $(\alpha ,d,\beta )$-su\-per\-pro\-cess $X$ have the same
distribution. Now, for any $K>0$, define the stopping time 
\begin{equation}
\tau _{K}=\inf \{t>0:\left\vert \Delta X_{t}\right\vert >K\}.
\label{defstopintime}
\end{equation}%
In Section~\ref{sec:lj} we will show that if we define the process which
evolves as $X$ up to time $\tau _{K}$ and then continues to evolve as $Y^{K}$
starting at $X_{\tau _{K}-}$, then this process has the same law as $Y^{K}$.
This together with the fact that $\tau _{K}\uparrow \infty $ as $%
K\rightarrow \infty $ (see Lemma~\ref{convtau}) implies that it is enough to
show existence of the SILT for the process $Y^{K}$. This task will be
accomplished in Section~\ref{sec:Y_SILT}, modulo some technical moment
estimates that will be derived in Section~\ref{moment}. The main steps
leading to the proof of Theorem~\ref{Theomain} will be described in the next
section.

\section{Proof of Theorem \protect\ref{Theomain}}

\label{sec:lj} The process $Y^{K}$ whose Laplace transform is given by (\ref%
{equt:3}), (\ref{equat1pp}) has the following decomposition: 
\begin{eqnarray}
Y_{t}^{K}(\varphi ) &=&\mu (\varphi )+\int_{0}^{t}Y_{s}^{K}(\Delta _{\alpha
}\varphi )ds-C_{\beta }(K)\int_{0}^{t}Y_{s}^{K}(\varphi )ds  \notag \\
&&+\int_{0}^{t}\int_{0}^{K}\int_{\mathbb{R}^{d}}r\varphi (x)\widetilde{N}%
_{Y^{K}}(dx,dr,ds),\;\;\forall t\geq 0,  \label{descompparayk}
\end{eqnarray}%
where $\widetilde{N}_{Y^{K}}=N_{Y^{K}}-\hat{N}_{Y^{K}}$, and 
\begin{eqnarray*}
N_{Y^{K}}(dx,dr,ds) &=&\sum_{\{(x,r,s):\Delta Y_{s}^{K}=r\delta
_{x}\}}\delta _{(x,r,s)}, \\
\hat{N}_{Y^{K}}(dx,dr,ds) &=&\eta 1_{(0,K]}(r)r^{-2-\beta }drY_{s}^{K}(dx)ds.
\end{eqnarray*}%
Note that $N_{Y^{K}}$ is defined in a way analogous to (\ref{dmc}), however
it does not have jumps \textquotedblleft greater\textquotedblright\ than $K$.

In the following lemma we are going to construct the probability space where 
$Y^K$ coincides with $X$ up to the stopping time $\tau_K$.

\begin{lemma}
\label{samedisyandyp} There exists a probability space on which a pair of
processes $(\ddot{Y}^K,X)$ is defined and possesses the following properties:

\begin{itemize}
\item[\textbf{(a)}] $\ddot{Y}^K$ coincides in law with $Y^K$,

\item[\textbf{(b)}] $\ddot{Y}^K_t=X_t,\;\;\;\forall t<\tau_K\,.$
\end{itemize}
\end{lemma}

\begin{proof}
Define 
\begin{equation*}
\Omega \equiv \Omega ^{\prime }\times \Omega ^{\circ },\;\;\mathcal{F}\equiv 
\mathcal{F}^{\prime }\times \mathcal{F}^{\circ },\;\;\mathcal{F}_{t}\equiv 
\mathcal{F}_{t}^{\prime }\times \mathcal{F}_{t}^{\circ },
\end{equation*}%
and let 
\begin{equation*}
\ddot{Y}_{t}^{K}(w^{\prime },w^{\circ })\equiv \left\{ 
\begin{tabular}{ll}
$X_{t}(w^{\prime })$, & $t<\tau _{K}(w^{\prime })$, \\ 
$w^{\circ }(t-\tau _{K}(w^{\prime }))$, & $t\geq \tau _{K}(w^{\prime }).$%
\end{tabular}%
\ \right.
\end{equation*}%
Define the measure $P$ on $(\Omega ,\mathcal{F})$: 
\begin{equation*}
P(B\times C)=\int_{\Omega ^{\prime }}1_{B}(w^{\prime })P^{\tau
_{K}(w^{\prime })}(C)P^{\prime }(dw^{\prime }),
\end{equation*}%
where%
\begin{equation}
P^{\tau _{K}(w^{\prime })}(C)=Q_{X_{\tau _{K}(\omega ^{\prime
})-}}(\{w^{\circ }\in \Omega ^{\circ }:\ddot{Y}^{K}\left( w^{\prime
},w^{\circ })\in C\}\right) .  \label{defofpt}
\end{equation}

Let $\varphi \in $\ {Dom} $(\Delta _{\alpha })$ and $t>0$. From the
definition of $\ddot{Y}^{K}$ we have%
\begin{eqnarray}
\ddot{Y}_{t}^{K}(\varphi ) &=&\mu (\varphi )+\int_{0}^{t}\ddot{Y}%
_{s}^{K}(\Delta _{\alpha }\varphi )ds-K_{\beta }\int_{0}^{t}\ddot{Y}%
_{s}^{K}(\varphi )ds  \notag \\
&&+\int_{0}^{t}\int_{0}^{K}\int_{\mathbb{R}^{d}}r\varphi (x)\widetilde{N}_{%
\ddot{Y}^{K}}(dx,dr,ds),  \label{mypp}
\end{eqnarray}%
where $N_{\ddot{Y}^{K}}$ is defined by (\ref{dmc}) for $t<\tau _{K}$ and 
\begin{equation*}
N_{\ddot{Y}^{K}}(dx,dr,ds)=\sum_{\{(x,r,s):\Delta \ddot{Y}_{s}^{K}=r\delta
_{x}\}}\delta _{(x,r,s)},
\end{equation*}%
for $t\geq \tau _{K}$. Let us check that $\int_{0}^{t}\int_{0}^{K}\int_{%
\mathbb{R}^{d}}r\varphi (x)\widetilde{N}_{\ddot{Y}^{K}}(dx,dr,ds)$ is an ${%
\mathcal{\ F}}_{t}$-martingale: For any $t>u$, $B\in \mathcal{F}_{u}^{\prime
}$, $C\in \mathcal{F}_{u}^{\circ }$, we obtain by using the definition (\ref%
{defofpt}) of $P^{\tau _{K}(w^{\prime })}$%
\begin{align*}
& P\left( 1_{B\times C}\left( \int_{0}^{t}\int_{0}^{K}\int_{\mathbb{R}%
^{d}}r\varphi (x)\widetilde{N}_{\ddot{Y}^{K}}(dx,dr,ds)\right. \right. \\
& \left. \left. -\int_{0}^{u}\int_{0}^{K}\int_{\mathbb{R}^{d}}r\varphi (x)%
\widetilde{N}_{\ddot{Y}^{K}}(dx,dr,ds)\right) \right) \\
& =\int_{B}P^{\tau _{K}(w^{\prime })}\left( 1_{C}\left(
\int_{0}^{t}\int_{0}^{K}\int_{\mathbb{R}^{d}}r\varphi (x)\widetilde{N}_{%
\ddot{Y}^{K}}(dx,dr,ds)\right. \right. \\
& \left. \left. -\int_{0}^{u}\int_{0}^{K}\int_{\mathbb{R}^{d}}r\varphi (x)%
\widetilde{N}_{\ddot{Y}^{K}}(dx,dr,ds)\right) \right) P^{\prime }(dw^{\prime
}) \\
& =\int_{B}P^{\tau _{K}(w^{\prime })}\left( 1_{C}\left(
\int_{0}^{t}\int_{0}^{K}\int_{\mathbb{R}^{d}}r\varphi (x)\widetilde{N}_{%
\ddot{Y}^{K}}(dx,dr,ds)\right. \right. \\
& \left. \left. -\int_{0}^{(t\wedge \tau _{K}(w^{\prime }))\vee
u}\int_{0}^{K}\int_{\mathbb{R}^{d}}r\varphi (x)\widetilde{N}_{\ddot{Y}%
^{K}}(dx,dr,ds)\right) \right) P^{\prime }(dw^{\prime }) \\
& +\int_{B}P^{\tau _{K}(w^{\prime })}\left( 1_{C}\left( \int_{0}^{(t\wedge
\tau _{K}(w^{\prime }))\vee u}\int_{0}^{K}\int_{\mathbb{R}^{d}}r\varphi (x)%
\widetilde{N}_{\ddot{Y}^{K}}(dx,dr,ds)\right. \right.
\end{align*}%
\begin{eqnarray*}
&&\left. \left. -\int_{0}^{u}\int_{0}^{K}\int_{\mathbb{R}^{d}}r\varphi (x)%
\widetilde{N}_{\ddot{Y}^{K}}(dx,dr,ds)\right) \right) P^{\prime }(dw^{\prime
}) \\
&=&\int_{B}P^{\tau _{K}(w^{\prime })}\left( Q_{X_{\tau _{K}(w^{\prime
})-}}[1_{C}\right. \\
&&\mbox{}\hspace*{1cm}\times \left. \int_{(t\wedge \tau _{K}(w^{\prime
}))\vee u}^{t}\int_{0}^{K}\int_{\mathbb{R}^{d}}r\varphi (x)\widetilde{N}_{%
\ddot{Y}^{K}}(dx,dr,ds)|\mathcal{F}_{t\wedge \tau _{K}(w^{\prime }))\vee
u}^{\circ }]\right) P^{\prime }(dw^{\prime }) \\
&&+\int_{B}P^{\tau _{K}(w^{\prime })}\left( 1_{C}\right. \\
&&\mbox{}\hspace*{1cm}\left. \times \int_{u}^{(t\wedge \tau _{K}(w^{\prime
}))\vee u}\int_{0}^{K-}\int_{\mathbb{R}^{d}}r\varphi (x)\widetilde{N}%
_{X}(dx,dr,ds)\right) P^{\prime }(dw^{\prime }) \\
&=&\int_{B}P^{\tau _{K}(w^{\prime })}\left( 1_{C}Q_{X_{\tau _{K}(w^{\prime
})-}}\left( \int_{(t\wedge \tau _{K}(w^{\prime }))\vee
u}^{t}\int_{0}^{K}\int_{\mathbb{R}^{d}}r\varphi (x)\widetilde{N}_{\ddot{Y}%
^{K}}(dx,dr,ds)\right. \right. \\
&&\mbox{}\hspace*{1cm}\left. \left. \left\vert \mathcal{F}_{(t\wedge \tau
_{K}(w^{\prime }))\vee u}^{\circ }\right. \right) \right) P^{\prime
}(dw^{\prime }) \\
&&+\int_{B}P^{\tau _{K}(w^{\prime })}(C) \\
&&\mbox{}\hspace*{1cm}\times \int_{u}^{(t\wedge \tau _{K}(w^{\prime }))\vee
u}\int_{0}^{K-}\int_{\mathbb{R}^{d}}r\varphi (x)\widetilde{N}%
_{X}(dx,dr,ds)P^{\prime }(dw^{\prime }) \\
&=&\int_{B}1_{C}(X_{\cdot \wedge u}(\omega ^{\prime }))1_{\{u<\tau _{K}\}} \\
&&\mbox{}\hspace*{1cm}\times P^{\prime }\left( \int_{u}^{(t\wedge \tau
_{K}(w^{\prime }))\vee u}\int_{0}^{K-}\int_{\mathbb{R}^{d}}r\varphi (x)%
\widetilde{N}_{X}(dx,dr,ds)|\mathcal{F}_{u}^{\prime }\right) P^{\prime
}(dw^{\prime }),
\end{eqnarray*}%
where in the last equality for the first term we have used the fact that for %
\mbox{$P'$-a.s $\omega'$}, $\tilde{N}_{\ddot{Y}^{K}}$ is a $(Q_{X_{\tau
_{K}(w^{\prime })-}},\mathcal{F}_{t})$-martingale measure on $\mathbb{R}%
^{d}\times \mathbb{R}_{+}\times \lbrack \tau _{K}\wedge t,t]$. 
As for the second term we have used the simple identity 
\begin{equation*}
P^{\tau _{K}(w^{\prime })}(C)1_{\{u<\tau _{K}\}}=1_{C}(X_{\cdot \wedge
u}(\omega ^{\prime }))1_{\{u<\tau _{K}\}}
\end{equation*}%
for any $C\in \mathcal{F}_{u}^{\circ }$. Now use the fact that $\widetilde{N}%
_{X}$ is a $(P^{\prime },\mathcal{F}_{t}^{\prime })$-martingale measure to
get that 
\begin{equation*}
P^{\prime }\left( \int_{u}^{(t\wedge \tau _{K}(w^{\prime }))\vee
u}\int_{0}^{K-}\int_{\mathbb{R}^{d}}r\varphi (x)\widetilde{N}_{X}(dx,dr,ds)|%
\mathcal{F}_{u}^{\circ }\right) =0,
\end{equation*}%
and the proof that $\int_{0}^{t}\int_{0}^{K}\int_{\mathbb{R}^{d}}r\varphi (x)%
\widetilde{N}_{\ddot{Y}^{K}}(dx,dr,ds)$ is a martingale is complete.

Then, due to the uniqueness of the decomposition (\cite{K-R}, Theorem 7) we
conclude that $\ddot{Y}^{K}$ has the same distribution as $Y^{K}$.
\end{proof}

\paragraph{\textbf{Convention}}

Based on the above lemma, from now on we will assume that $Y^K, X$ are
defined on the same probability space and $Y^K_t=X_t\,, \forall t<\tau_K\,.$

\bigskip

Now we are going to show that time $\tau_K$ can be made greater than any
constant $T$ with probability arbitrary close to $1$ by taking $K$
sufficiently large.

\begin{lemma}
\label{convtau} For every $T>0$ and $\varepsilon >0$, there exists $K>0$
such that \mbox{$P(\tau _{K}\leq T)\leq \varepsilon .$}
\end{lemma}

\begin{proof}
Let $Z_{T}^{K}$ the number of jumps of height greater than $K$ in $%
[0,T]\times \mathbb{R}^{d}$, that is $Z_{T}^{K}=N([K,+\infty )\times \lbrack
0,T]\times \mathbb{R}^{d})$. Then there exists (see \cite{M-P}, page 1430) a
standard Poisson process $A_{t}^{K}$ such that 
\begin{equation*}
Z_{T}^{K}=A_{c_{\beta }K^{-1-\beta }\int_{0}^{T}X_{s}(\mathbb{R}^{d})ds}^{K},
\end{equation*}%
for some positive constant $c_{\beta }$. Then from the Markov inequality we
have, 
\begin{eqnarray*}
P(\tau _{K}\leq T) &=&P(Z_{T}^{K}\geq 1) \\
&=&P\left( A_{c_{\beta }K^{-1-\beta }\int_{0}^{T}X_{s}(\mathbb{R}%
^{d})ds}^{K}\geq 1\right) \\
&=&P\left( A_{c_{\beta }K^{-1-\beta }\int_{0}^{T}X_{s}(\mathbb{R}%
^{d})ds}^{K}\geq 1,\frac{c_{\beta }}{K^{1+\beta }}\int_{0}^{T}X_{s}(\mathbb{R%
}^{d})ds\geq K^{-1}\right) \\
&&+P\left( A_{c_{\beta }K^{-1-\beta }\int_{0}^{T}X_{s}(\mathbb{R}%
^{d})ds}^{K}\geq 1,\frac{c_{\beta }}{K^{1+\beta }}\int_{0}^{T}X_{s}(\mathbb{R%
}^{d})ds<K^{-1}\right) \\
&\leq &P\left( \int_{0}^{T}X_{s}(\mathbb{R}^{d})\geq c_{\beta }^{-1}K^{\beta
}\right) +P\left( A_{K^{-1}}^{K}\geq 1\right) \\
&\leq &c_{\beta }K^{-\beta }E_{\mu }\left[ \int_{0}^{T}X_{s}(\mathbb{R}^{d})%
\right] +\left( 1-P\left( A_{K^{-1}}^{K}=0\right) \right) \\
&=&c_{\beta }T\mu (\mathbb{R}^{d})K^{-\beta }+\left( 1-\exp \left(
-K^{-1}\right) \right) .
\end{eqnarray*}%
The result follows, since the right hand side goes to $0$ as $K\rightarrow
\infty $.
\end{proof}

Now the proof of Theorem \ref{Theomain} relies on the following proposition.

\begin{proposition}
\label{TeoconL1}Let $K>0$ and $Y^{K}$ be the truncated $(\alpha ,d,\beta )$%
-superprocess with initial measure $Y_{0}^{K}(dx)=\mu (dx)=h(x)dx$, where $h$
is bounded and integrable with respect to Lebesgue measure $dx$ on $\mathbb{R%
}^{d}$.

\begin{itemize}
\item[\textbf{(a)}] For $d/2<\alpha $ there exists a process $\gamma
_{Y^{K}}^{K}=\{\gamma _{Y^{K}}^{K}(T):T\geq 0\}$ such that 
\begin{equation*}
\lim_{\varepsilon \downarrow 0}E\left[ \sup_{t<T}\left\vert \gamma
_{Y^{K},\varepsilon }^{K}(t))-\gamma _{Y^{K}}^{K}(t))\right\vert \right]
=0,\;\;\forall T>0,
\end{equation*}%
and for any $\lambda >0,$ 
\begin{eqnarray*}
\gamma _{Y^{K}}^{K}(T) &=&\lambda \int_{0}^{T}\int_{0}^{t}\int_{\mathbb{R}%
^{2d}}G^{\lambda }(x-y)Y_{s}^{K}(dx)Y_{t}^{K}(dy)dtds \\
&&-\int_{0}^{T}\int_{\mathbb{R}^{d}}Y_{T}^{K}(G^{\lambda }(x-\cdot
))Y_{s}^{K}(dx) \\
&&+\int_{0}^{T}\int_{\mathbb{R}^{d}}Y_{s}^{K}(G^{\lambda }(x-\cdot
))Y_{s}^{K}(dx)ds \\
&&+\int_{0}^{T}\int_{\mathbb{R}^{d}}\int_{0}^{t}\int_{\mathbb{R}%
^{d}}G^{\lambda }(x-y)M^{K}(ds,dy)Y_{t}^{K}(dx)dt,\;\;\mathrm{a.s.}
\end{eqnarray*}

\item[\textbf{(b)}] For $d/(2+(1+\beta )^{-1})<\alpha \leq d/2$ there exists
a process $\tilde{\gamma}_{Y^{K}}^{K}=\{\tilde{\gamma}_{Y^{K}}^{K}(T):T\geq
0\}$ such that 
\begin{equation*}
\lim_{\varepsilon \downarrow 0}E\left[ \sup_{t<T}\left\vert \tilde{\gamma}%
_{Y^{K},\varepsilon }^{K}(t))-\tilde{\gamma}_{Y^{K}}^{K}(t))\right\vert %
\right] =0,\;\;\forall T>0,
\end{equation*}%
and for any $\lambda >0,$ 
\begin{eqnarray*}
\tilde{\gamma}_{Y^{K}}^{K}(T) &=&\lambda \int_{0}^{T}\int_{0}^{t}\int_{%
\mathbb{R}^{2d}}G^{\lambda }(x-y)Y_{s}^{K}(dx)Y_{t}^{K}(dy)dtds \\
&&-\int_{0}^{T}\int_{\mathbb{R}^{d}}Y_{T}^{K}(G^{\lambda }(x-\cdot
))Y_{s}^{K}(dx)ds \\
&&+\int_{0}^{T}\int_{\mathbb{R}^{d}}\int_{0}^{t}\int_{\mathbb{R}%
^{d}}G^{\lambda }(x-y)M^{K}(ds,dy)Y_{t}^{K}(dx)dt,\;\;\mathrm{a.s.}
\end{eqnarray*}
\end{itemize}
\end{proposition}

\paragraph{Proof}

Postponed.

\mbox{}\newline
The above proposition immediately yields:

\begin{proof}[Proof of Theorem \protect\ref{Theomain}]
Fix arbitrary $\varepsilon,\delta>0$ and let $d/2<\alpha$. Since $X_t=Y^K_t$
for any $t<\tau_K$, we immediately get that 
\begin{equation*}
\gamma _{X}(t)=\gamma _{Y^K}(t), \;\;\forall t<\tau_K\,,
\end{equation*}
and $\gamma _{X}(t)$ satisfies Tanaka formula~(\ref{tkfmala}) for $t<\tau_K$%
. Moreover, since by Lemma~\ref{convtau}, $\tau_K\uparrow \infty$, as $%
K\rightarrow\infty$, there is no problem to define $\gamma _{X}(t)$
satisfying~(\ref{tkfmala}) for any $t>0$.

Now let us check the convergence part of the theorem. For any $T>0$, by
Lemma \ref{convtau}, we can fix $K>0$ such that $P(\tau _{K}\leq T)\leq
\delta $. Then 
\begin{eqnarray*}
\lefteqn{\lim_{\varepsilon _{1}\downarrow 0}P\left( \sup_{t\leq T}\left\vert
\gamma _{X,\varepsilon _{1}}(t)-\gamma _{X}(t)\right\vert >\varepsilon
\right) } \\
&\leq &\lim_{\varepsilon _{1}\downarrow 0}P\left( \sup_{t\leq T}\left\vert
\gamma _{X,\varepsilon _{1}}(t)-\gamma _{X}(t)\right\vert >\varepsilon ,\tau
_{K}>T\right) +P\left( \tau _{K}\leq T\right) \\
&\leq &\lim_{\varepsilon _{1}\downarrow 0}P\left( \sup_{t\leq T}\left\vert
\gamma _{Y^{K},\varepsilon _{1}}(t)-\gamma _{Y^{K}}(t)\right\vert
>\varepsilon ,\tau _{K}>T\right) +\delta \\
&=&\delta ,
\end{eqnarray*}%
and since $\delta ,\varepsilon >0$ were arbitrary the proof of convergence
is complete.

The proof of part~\textbf{(b)} of the theorem goes along the same lines.
\end{proof}

\section{Existence of SILT for $Y^K$ --- proof of Proposition~\protect\ref%
{TeoconL1}}

\label{sec:Y_SILT} 
Fix arbitrary $K>0$. First, we derive very useful moment estimates for $Y^K$%
. Let $\{S_{t}^{K}:t\geq 0\}$ denote the solution of the partial
differential equation 
\begin{equation*}
\frac{\partial v_{t}^{K}}{\partial t}=(\Delta _{\alpha }-C_{\beta
}(K))v_{t}^{K}.
\end{equation*}
That is, $\{S_{t}^{K}:t\geq 0\}$ is the semigroup defined as 
\begin{equation}
S_{t}^{K}=e^{-C_{\beta }(K)t}S_{t}.  \label{relstrunsem}
\end{equation}
Notice that $S_{t}^{K}\varphi \leq S_{t}\varphi $, for all non-negative
bounded measurable functions $\varphi $.

Following Theorem 3.1 of \cite{I} we have that for $\varphi ,\psi \in 
\mathcal{B}_{b}(\mathbb{R}^{d}),$%
\begin{equation}
E\left[ \left. \exp \left( -Y_{t}^{K}(\varphi )-\int_{0}^{t}Y_{s}^{K}(\psi
)ds\right) \right\vert Y_{0}^{K}=\mu \right] =\exp \left( -\mu
(V_{t}^{K}(\varphi ,\psi ))\right) ,  \label{lapfuntrun}
\end{equation}%
where $\mu \in \mathcal{M}_{F}(\mathbb{R}^{d})$ and $V_{t}^{K}(\varphi ,\psi
)$ denotes the unique non-negative solution to the following evolution
equation 
\begin{equation}
v_{t}^{K}=S_{t}^{K}\varphi +\int_{0}^{t}S_{s}^{K}(\psi
)ds-\int_{0}^{t}S_{t-s}^{K}\left( \Phi ^{K}(v_{s}^{K})\right) ds,\text{ \ \ }%
t\geq 0,  \label{evoltrunc}
\end{equation}%
where 
\begin{equation}
\Phi ^{K}(x)=\sum_{m=2}^{\infty }\frac{(-1)^{m}}{m!}\chi (m)x^{m}
\label{seriesdeftrunfi}
\end{equation}%
and 
\begin{equation}
\chi (m)=\frac{\eta K^{m-1-\beta }}{m-1-\beta }.  \label{equt:12}
\end{equation}

Now we are going to calculate the first two moments of $Y^{K}$.

\begin{lemma}
\label{lemamome}Let $\varphi $ be a non-negative function on $\mathcal{B}%
_{b}(\mathbb{R}^{d})$ and $t>0$. Then 
\begin{eqnarray}
E_{\mu }\left[ Y_{t}^{K}(\varphi )\right] &=&\mu \left( S_{t}^{K}\varphi
\right) ,  \label{smfmla} \\
E_{\mu }\left[ \left( Y_{t}^{K}(\varphi )\right) ^{2}\right] &=&\left( \mu
\left( S_{t}^{K}\varphi \right) \right) ^{2}+\chi (2)\mu \left(
\int_{0}^{t}S_{t-s}^{K}\left( \left( S_{s}^{K}\varphi \right) ^{2}\right)
ds\right) .  \notag
\end{eqnarray}
\end{lemma}

\begin{proof}
From (\ref{lapfuntrun}) we have 
\begin{equation}
E_{\mu }\left[ e^{-\lambda Y_{t}^{K}(\varphi )}\right] =e^{-\mu \left(
v_{t}^{K}(\lambda )\right) }  \label{laptrunc}
\end{equation}%
where 
\begin{equation}
v_{t}^{K}(\lambda )=\lambda S_{t}^{K}\varphi -\int_{0}^{t}S_{t-s}^{K}\left(
\Phi ^{K}(v_{s}^{K}(\lambda ))\right) ds.  \label{integeqc2}
\end{equation}%
Using the elementary inequality $e^{-x}-1+x\leq x^{2}/2$, $x\geq 0$, and (%
\ref{defapxintpsim}) we have%
\begin{equation*}
\Phi ^{K}(x)\leq \frac{\chi (2)}{2}x^{2},\ \ x\geq 0.
\end{equation*}%
Let $||\cdot ||_{\infty }$ be the supremum norm, then $0\leq
v_{t}^{K}(\lambda )\leq \lambda S_{t}^{K}\varphi \leq \lambda ||\varphi
||_{\infty }$ and the previous inequality implies 
\begin{equation*}
\int_{0}^{t}S_{t-s}^{K}\left( \Phi ^{K}(v_{s}^{K}(\lambda ))\right) ds\leq 
\frac{\chi (2)||\varphi ||_{\infty }^{2}t}{2}\times \lambda ^{2}.
\end{equation*}%
Further from (\ref{integeqc2}) we get 
\begin{eqnarray*}
\lim_{\lambda \downarrow 0}\frac{\lambda S_{t}^{K}\varphi -v_{t}^{K}(\lambda
)}{\lambda } &=&\lim_{\lambda \downarrow 0}\frac{1}{\lambda }%
\int_{0}^{t}S_{t-s}^{K}\left( \Phi ^{K}(v_{s}^{K}(\lambda ))\right) ds \\
&\leq &\lim_{\lambda \downarrow 0}\frac{\chi (2)||\varphi ||_{\infty }^{2}t}{%
2}\times \lambda =0,
\end{eqnarray*}%
and we write this like%
\begin{equation}
v_{t}^{K}(\lambda )=\lambda S_{t}^{K}\varphi -o(\lambda ).  \label{extsol1}
\end{equation}%
This implies 
\begin{eqnarray*}
E_{\mu }\left[ Y_{t}^{K}(\varphi )\right] &=&\lim_{\lambda \downarrow 0}%
\frac{1-E_{\mu }\left[ e^{-\lambda Y_{t}^{K}(\varphi )}\right] }{\lambda } \\
&=&\lim_{\lambda \downarrow 0}\frac{1-e^{-\lambda \mu \left(
S_{t}^{K}\varphi \right) +o(\lambda )}}{\lambda } \\
&=&\lim_{\lambda \downarrow 0}\frac{1-e^{-\lambda \mu \left(
S_{t}^{K}\varphi \right) +o(\lambda )}}{\lambda \mu \left( S_{t}^{K}\varphi
\right) -o(\lambda )}\lim_{\lambda \downarrow 0}\frac{\lambda \mu \left(
S_{t}^{K}\varphi \right) -o(\lambda )}{\lambda }=\mu \left( S_{t}^{K}\varphi
\right) .
\end{eqnarray*}%
Now, to calculate the second moment we follow the ideas used in the proof of
Proposition 11 of Chapter II from~\cite{L}:%
\begin{eqnarray*}
E_{\mu }\left[ \left( Y_{t}^{K}(\varphi )\right) ^{2}\right]
&=&\lim_{\lambda \downarrow 0}\frac{2}{\lambda ^{2}}E_{\mu }\left[
e^{-\lambda Y_{t}^{K}(\varphi )}-1+\lambda Y_{t}^{K}(\varphi )\right] \\
&=&\lim_{\lambda \downarrow 0}\frac{2}{\lambda ^{2}}\left( e^{-\mu
(v_{t}^{K}(\lambda ))}-1+\lambda \mu \left( S_{t}^{K}\varphi \right) \right)
\\
&=&\lim_{\lambda \downarrow 0}\frac{2}{\lambda ^{2}}\left( e^{-\lambda \mu
\left( S_{t}^{K}\varphi \right) +\mu \left( \int_{0}^{t}S_{t-s}^{K}\left(
\Phi ^{K}(v_{s}^{K}(\lambda ))\right) ds\right) }-1+\lambda \mu \left(
S_{t}^{K}\varphi \right) \right) \\
&=&\lim_{\lambda \downarrow 0}\frac{2}{\lambda ^{2}}(\sum_{n=0}^{\infty }%
\frac{1}{n!}\left( \mu \left( \int_{0}^{t}S_{t-s}^{K}\left( \Phi
^{K}(v_{s}^{K}(\lambda ))\right) ds\right) -\lambda \mu \left(
S_{t}^{K}\varphi \right) \right) ^{n} \\
&&-1+\lambda \mu \left( S_{t}^{K}\varphi \right) ).
\end{eqnarray*}

Using the series expansion (\ref{seriesdeftrunfi}) for $\Phi ^{K}$ and (\ref%
{extsol1}) we obtain%
\begin{eqnarray*}
\int_{0}^{t}S_{t-s}^{K}\left( \Phi ^{K}(v_{s}^{K}(\lambda ))\right) ds
&=&\int_{0}^{t}S_{t-s}^{K}(\Phi ^{K}(\lambda S_{s}^{K}\varphi -o(\lambda
)))ds \\
&=&\int_{0}^{t}S_{t-s}^{K}\left( \frac{\chi (2)}{2!}\lambda ^{2}\left(
S_{s}^{K}\varphi \right) ^{2}+o(\lambda ^{2})\right) ds \\
&=&\frac{\chi (2)}{2}\lambda ^{2}\int_{0}^{t}S_{t-s}^{K}\left( \left(
S_{s}^{K}\varphi \right) ^{2}\right) ds+o(\lambda ^{2}).
\end{eqnarray*}

Then 
\begin{align*}
& E_{\mu }\left[ \left( Y_{t}^{K}(\varphi )\right) ^{2}\right] \\
& =\lim_{\lambda \downarrow 0}\frac{2}{\lambda ^{2}}\left( \frac{\chi (2)}{2}%
\lambda ^{2}\mu \left( \int_{0}^{t}S_{t-s}^{K}\left( \left( S_{s}^{K}\varphi
\right) ^{2}\right) ds\right) +o(\lambda ^{2})\right. \\
& \left. +\frac{1}{2!}\left( \frac{\chi (2)}{2}\lambda ^{2}\mu \left(
\int_{0}^{t}S_{t-s}^{K}\left( \left( S_{s}^{K}\varphi \right) ^{2}\right)
ds\right) -\lambda \mu \left( S_{t}^{K}\varphi \right) +o(\lambda
^{2})\right) ^{2}\right) \\
& =\lim_{\lambda \downarrow 0}\frac{2}{\lambda ^{2}}\left( \frac{\chi (2)}{2}%
\lambda ^{2}\mu \left( \int_{0}^{t}S_{t-s}^{K}\left( \left( S_{s}^{K}\varphi
\right) ^{2}\right) ds\right) +\frac{1}{2}\lambda ^{2}\left( \mu \left(
S_{t}^{K}\varphi \right) \right) ^{2}+o(\lambda ^{2})\right) \\
& =\chi (2)\mu \left( \int_{0}^{t}S_{t-s}^{K}\left( \left( S_{s}^{K}\varphi
\right) ^{2}\right) ds\right) +\left( \mu \left( S_{t}^{K}\varphi \right)
\right) ^{2},
\end{align*}%
and we are done.
\end{proof}

\begin{remark}
Using binary directed graphs, Dynkin in \cite{D} gives a formula for the
moments of supeprocesses, where the Laplace functional (\ref{lapfuntrun})
has an evolution equation (\ref{evoltrunc}) with only one term $m=2$ in (\ref%
{seriesdeftrunfi}). For the $Y^{K}$ superprocess it is also possible, but
here the main difference is that the directed graphs are not necessarily
binary.
\end{remark}

\begin{corollary}
\label{cor:1} Let $\varphi ,\psi $ be non-negative functions on $\mathcal{B}%
_{b}(\mathbb{R}^{d})$ and $t\geq s>0$. Then 
\begin{eqnarray*}
E_{\mu }\left[ Y_{t}^{K}(\varphi )Y_{s}^{K}(\psi )\right] &=&\mu \left(
S_{t}^{K}\varphi \right) \mu \left( S_{s}^{K}\psi \right) \\
&&+\chi (2)\mu \left( \int_{0}^{s}S_{r}^{K}(S_{t-r}^{K}\varphi
S_{s-r}^{K}\psi )dr\right) .
\end{eqnarray*}
\end{corollary}

\begin{proof}
First, use the Markov property for $Y^{K}$ to get 
\begin{equation*}
E_{\mu }\left[ Y_{t}^{K}(\varphi )|\mathcal{F}_{s}\right] =Y_{s}^{K}\left(
S_{t-s}^{K}\varphi \right) .
\end{equation*}%
Therefore,%
\begin{align*}
& E_{\mu }\left[ Y_{t}^{K}(\varphi )Y_{s}^{K}(\psi )\right] \\
& =E_{\mu }\left[ Y_{s}^{K}\left( S_{t-s}^{K}\varphi \right) Y_{s}^{K}(\psi )%
\right] \\
& =\frac{1}{4}\left( \left( E_{\mu }\left[ Y_{s}^{K}\left(
S_{t-s}^{K}\varphi \right) +Y_{s}^{K}(\psi )\right] \right) ^{2}-\left(
E_{\mu }\left[ Y_{s}^{K}\left( S_{t-s}^{K}\varphi \right) -Y_{s}^{K}(\psi )%
\right] \right) ^{2}\right) ,
\end{align*}

and we are done by Lemma \ref{lemamome}.
\end{proof}

\begin{corollary}
\label{cor:2} Let $\varphi $ be non-negative functions on $\mathcal{B}_{b}(%
\mathbb{R}^{d}\times \mathbb{R}^{d})$ and $t\geq s>0$. Then 
\begin{eqnarray*}
\lefteqn{E_{\mu }\left[ \int_{\mathbb{R}^{d}\times \mathbb{R}^{d}}\varphi
(x,y)Y_{t}^{K}(dx)Y_{s}^{K}(dy)\right] } \\
&=&\int_{\mathbb{R}^{4d}}\mu (dx_{1})\mu
(dx_{2})p_{t}(z_{1}-x_{1})p_{s}(z_{2}-x_{2})\varphi (z_{1},z_{2})dz_{1}dz_{2}
\\
&&+\chi (2)\int_{0}^{s}\int_{\mathbb{R}^{4d}}\mu
(dx)p_{r}(y-x)dyp_{t-r}(z_{1}-y)p_{s-r}(z_{2}-y)\varphi
(z_{1},z_{2})dz_{1}dz_{2}dr.
\end{eqnarray*}
\end{corollary}

\begin{proof}
Use Corollary~\ref{cor:1} and approximation of the $\psi (x,y)$ by functions
in the form $\sum_{i}\varphi _{i}(x)\phi _{i}(y)$ to derive the result. We
leave the details to the reader.
\end{proof}

Next proposition gives bounds on some fractional moments of $Y^K$ and
requires much more work than we have done in Lemma~\ref{lemamome}. Hence its
proof will be postponed till Section~\ref{moment}.

\begin{proposition}
\label{markcorolory}Let $1+\beta <p<2$ and $0<\varepsilon \leq 1$. If 
\begin{equation*}
d<\alpha \left( 2+\frac{1}{p}\right) ,
\end{equation*}%
then there exists a constant $c=c(K,p,d,\alpha ,\beta )$ such that 
\begin{equation*}
E_{\mu }\left[ \int_{\mathbb{R}^{d}}Y_{t}^{K}(p_{\varepsilon }(\cdot
-x))\left( \int_{0}^{t}Y_{s}^{K}(G^{\lambda }(\cdot \cdot -x))ds\right)
^{p}dx\right] <c(K,p,d,\alpha ,\beta ).
\end{equation*}%
Moreover 
\begin{equation}
E\left[ \int_{\mathbb{R}^{d}}\left( \int_{0}^{t}Y_{s}^{K}(G^{\lambda }(\cdot
\cdot -x))ds\right) ^{p}Y_{t}^{K}(dx)\right] <\infty .  \label{momreferee}
\end{equation}
\end{proposition}

\begin{remark}
When $K$ goes to infinity, then $\chi (2)$ goes to infinity, hence $%
c(K,p,d,\alpha ,\beta )$ goes to infinity and this is because $\chi (2)$ is
part of $c(K,p,d,\alpha ,\beta )$. The moment in (\ref{momreferee}) is
infinity when $K=+\infty $, because the $(\alpha ,d,\beta )$-superprocess
has moments of order less than $1+\beta $ and $p>1+\beta $.
\end{remark}

\paragraph{Proof}

Postponed. 

\mbox{}\newline
Now we can write the Tanaka-like formula for the approximating SILT of the
truncated superprocess $Y^{K}$. From Fubini theorem, (\ref{genagreeb}), (\ref%
{storepmtg}) and (\ref{mp}) (the martingale problem for the truncated
superprocess $Y^{K}$, \cite{K-R}) we have%
\begin{eqnarray}
\gamma _{Y^{K},\varepsilon }^{K}(T) &=&\int_{0}^{T}\int_{0}^{t}\int_{\mathbb{%
\ R}^{2d}}p_{\varepsilon }(x-y)Y_{s}^{K}(dx)Y_{t}^{K}(dy)dsdt  \notag \\
&=&\lambda e^{\lambda \varepsilon }\int_{0}^{T}\int_{0}^{t}\int_{\mathbb{R}%
^{2d}}G^{\lambda ,\varepsilon }(x-y)Y_{s}^{K}(dx)Y_{t}^{K}(dy)dtds  \notag \\
&&-e^{\lambda \varepsilon }\int_{0}^{T}\int_{\mathbb{R}^{d}}\int_{s}^{T}%
\int_{\mathbb{R}^{d}}\Delta _{\alpha }G^{\lambda ,\varepsilon
}(x-y)Y_{t}^{K}(dy)dtY_{s}^{K}(dx)ds  \notag \\
&=&\lambda e^{\lambda \varepsilon }\int_{0}^{T}\int_{0}^{t}\int_{\mathbb{R}%
^{2d}}G^{\lambda ,\varepsilon }(x-y)Y_{s}^{K}(dx)Y_{t}^{K}(dy)dsdt  \notag \\
&&-e^{\lambda \varepsilon }\int_{0}^{T}\int_{\mathbb{R}^{d}}Y_{T}^{K}(G^{%
\lambda ,\varepsilon }(x-\cdot ))Y_{s}^{K}(dx)ds  \notag \\
&&+e^{\lambda \varepsilon }\int_{0}^{T}\int_{\mathbb{R}^{d}}Y_{s}^{K}(G^{%
\lambda ,\varepsilon }(x-\cdot ))Y_{s}^{K}(dx)ds  \notag \\
&&+e^{\lambda \varepsilon }\int_{0}^{T}\int_{\mathbb{R}^{d}}\int_{s}^{T}%
\int_{\mathbb{R}^{d}}G^{\lambda ,\varepsilon
}(x-y)M^{K}(dt,dy)Y_{s}^{K}(dx)ds,  \label{takanaaprox}
\end{eqnarray}%
and 
\begin{eqnarray}
\tilde{\gamma}_{Y^{K},\varepsilon }^{K}(T) &=&\gamma _{Y^{K},\varepsilon
}^{K}(T)-e^{\lambda \varepsilon }\int_{0}^{T}\int_{\mathbb{R}%
^{d}}Y_{s}^{K}(G^{\lambda ,\varepsilon }(x-\cdot ))Y_{s}^{K}(dx)ds  \notag \\
&=&\lambda e^{\lambda \varepsilon }\int_{0}^{T}\int_{0}^{t}\int_{\mathbb{R}%
^{2d}}G^{\lambda ,\varepsilon }(x-y)Y_{s}^{K}(dx)Y_{t}^{K}(dy)dsdt  \notag \\
&&-e^{\lambda \varepsilon }\int_{0}^{T}\int_{\mathbb{R}^{d}}Y_{T}^{K}(G^{%
\lambda ,\varepsilon }(x-\cdot ))Y_{s}^{K}(dx)ds  \notag \\
&&+e^{\lambda \varepsilon }\int_{0}^{T}\int_{\mathbb{R}^{d}}\int_{s}^{T}%
\int_{\mathbb{R}^{d}}G^{\lambda ,\varepsilon
}(x-y)M^{K}(dt,dy)Y_{s}^{K}(dx)ds.  \label{takanaaprox1}
\end{eqnarray}%
Note that stochastic integrals ~in (\ref{takanaaprox}) and (\ref%
{takanaaprox1}) are well defined due to the moment bound given by
Proposition~\ref{markcorolory}. %
%
%
%
%

\begin{proof}[Proof of Proposition \protect\ref{TeoconL1}]
We are going to prove Proposition \ref{TeoconL1} via letting $\varepsilon
\rightarrow 0$ in~(\ref{takanaaprox}), and checking convergence of all the
terms. 
By Corollary~\ref{cor:2} and simple estimates we get%
\begin{eqnarray*}
&&E\left[ \int_{0}^{T}\int Y_{s}^{K}(G^{\lambda ,\varepsilon }(x-\cdot
))Y_{s}^{K}(dx)ds\right] \\
&\leq &\int_{0}^{T}\int_{\mathbb{R}^{4d}}\mu (dx_{1})\mu
(dx_{2})p_{s}(z_{1}-x_{1})p_{s}(z_{2}-x_{2})G^{\lambda
}(z_{1}-z_{2})dz_{1}dz_{2}ds \\
&&+\chi (2)\int_{0}^{T}\int_{0}^{s}\int_{\mathbb{R}^{4d}}\mu
(dx)p_{r}(y-x)dyp_{s-r}(z_{1}-y)p_{s-r}(z_{2}-y) \\
&&\times G^{\lambda }(z_{1}-z_{2})dz_{1}dz_{2}drds
\end{eqnarray*}%
\begin{align}
& \leq \int_{0}^{T}\int_{\mathbb{R}^{4d}}\mu (dx_{1})\mu
(dx_{2})\int_{0}^{\infty }e^{-\lambda
u}p_{u}(z_{1}-z_{2})p_{s}(z_{1}-x_{1})p_{s}(z_{2}-x_{2})dz_{1}dz_{2}duds 
\notag \\
& +\chi (2)\int_{0}^{T}\int_{0}^{s}\int_{\mathbb{R}^{4d}}\mu (dx)p_{r}(y-x)dy
\notag \\
& \times \int_{0}^{\infty }e^{-\lambda
u}p_{u}(z_{1}-z_{2})p_{s-r}(z_{1}-y)dz_{1}p_{s-r}(z_{2}-y)dz_{2}dudrds 
\notag \\
& \leq \int_{0}^{T}\int_{0}^{\infty }e^{-\lambda u}\int
p_{u+2s}(x_{1}-x_{2})\mu (dx_{1})\mu (dx_{2})duds  \notag \\
& +\chi (2)\int_{0}^{T}\int_{0}^{s}\int \mu (dx)p_{r}(y-x)dy\int_{0}^{\infty
}e^{-\lambda u}p_{u+2s-2r}(0)dudrds  \notag \\
& \leq \left\Vert h\right\Vert _{\infty }\mu (1)T\lambda ^{-1}+\mu (1)\chi
(2)\int_{0}^{T}\int_{0}^{s}\int_{0}^{\infty }e^{-\lambda
u}(u+2s-2r)^{-d/\alpha }dudrds,\;\;\;\forall T>0,  \label{calikf}
\end{align}%
where the last integral is convergent if $d<2\alpha $. Using~(\ref{calikf}),
the bound $G^{\lambda }\geq G^{\lambda ,\varepsilon }$ and the monotone
convergence theorem to get 
\begin{align}
& \lim_{\varepsilon \downarrow 0}E\left[ \sup_{t<T}\left\vert
\int_{0}^{t}\int Y_{s}^{K}(G^{\lambda }(x-\cdot )-G^{\lambda ,\varepsilon
}(x-\cdot ))Y_{s}^{K}(dx)ds\right\vert \right]  \notag \\
& \leq \lim_{\varepsilon \downarrow 0}E\left[ \int_{0}^{T}\int
Y_{s}^{K}(G^{\lambda }(x-\cdot )-G^{\lambda ,\varepsilon }(x-\cdot
))Y_{s}^{K}(dx)ds\right]  \notag \\
& =0,\;\;\;\forall T>0.  \label{leonosexqja}
\end{align}

In a similar way we can prove that 
\begin{equation}
\lim_{\varepsilon \downarrow 0}E\left[ \sup_{T<L}\left\vert
\int_{0}^{T}\int_{0}^{t}\int_{\mathbb{R}^{2d}}\left( G^{\lambda
}(x-y)-G^{\lambda ,\varepsilon }(x-y)\right)
Y_{s}^{K}(dx)Y_{t}^{K}(dy)dsdt\right\vert \right] =0,  \label{equt:5}
\end{equation}%
and 
\begin{equation}
\lim_{\varepsilon \downarrow 0}E\left[ \sup_{T<L}\left\vert
\int_{0}^{T}\int_{\mathbb{R}^{d}}Y_{T}^{K}(G^{\lambda }(x-\cdot )-G^{\lambda
,\varepsilon }(x-\cdot ))Y_{s}^{K}(dx)ds,\right\vert \right] =0,
\label{equt:6}
\end{equation}%
for all $L>0$ and $d<3\alpha $.

Now let us deal with the stochastic integral 
\begin{equation*}
\int_{0}^{T}\int F^{\varepsilon }(t,x)M^{K}(dt,dx),
\end{equation*}%
where 
\begin{equation*}
F^{\varepsilon }(t,x)=\int_{0}^{t}\int_{\mathbb{R}^{d}}G^{\lambda
,\varepsilon }(x-y)Y_{s}^{K}(dy)ds.
\end{equation*}%
This integral is well defined if (see~\cite{L-M}) 
\begin{equation}
E\left[ \left( \sum_{t\in J\cap \lbrack 0,T]}F^{\varepsilon }(t,\Delta
Y_{t}^{K})^{2}\right) ^{1/2}\right] <+\infty ,  \label{condintegr}
\end{equation}%
where $J$ denotes the set of all jump times of $X$. Let 
\begin{equation}
d<\alpha \left( 2+\frac{1}{1+\beta }\right) ,  \label{condplema}
\end{equation}%
hence we can choose $p\in (1+\beta ,2)$ such that 
\begin{equation}
d<\alpha \left( 2+\frac{1}{p}\right) .  \label{conauxp}
\end{equation}%
Since $p\in (1+\beta ,2)$ we can use the Jensen inequality to get 
\begin{equation*}
\left( \sum_{i\in I}a_{i}\right) ^{p/2}\leq \sum_{i\in I}a_{i}^{p/2},
\end{equation*}%
if $a_{i}\geq 0$ for all $i\in I$. This yields 
\begin{align*}
& E\left[ \left( \sum_{t\in J\cap \lbrack 0,T]}F^{\varepsilon }(t,\Delta
Y_{t}^{K})^{2}\right) ^{1/2}\right]  \\
& \leq \left( E\left[ \sum_{t\in J\cap \lbrack 0,T]}F^{\varepsilon
}(t,\Delta Y_{t}^{K})^{p}\right] \right) ^{1/p} \\
& =\left( E\left[ \eta \int_{0}^{T}\int \int_{0}^{K}u^{-\beta
-2}(F^{\varepsilon }(t,u\delta _{x}))^{p}duY_{t}^{K}(dx)dt\right] \right)
^{1/p} \\
& =c\left( \int_{0}^{T}E\left[ \int \left( \int_{0}^{t}\int G^{\lambda
,\varepsilon }(x-y)Y_{s}^{K}(dy)ds\right) ^{p}Y_{t}^{K}(dx)\right] dt\right)
^{1/p}.
\end{align*}%
Since $p$ satisfies~(\ref{conauxp}), the condition (\ref{condintegr})
follows from Proposition \ref{markcorolory}.

Let $F=F^{0}$. By Burkholder-Davis-Gundy inequality (see \cite{L-M}) and the
previous argument we get 
\begin{align}
& E\left[ \sup_{t\leq T}\left\vert \int_{0}^{t}\int (F(s,x)-F^{\varepsilon
}(s,x))M^{K}(ds,dx)\right\vert \right]  \notag \\
& \leq E\left[ \sup_{t\leq T}\int_{0}^{t}\int (\left\vert
F(s,x)-F^{\varepsilon }(s,x)\right\vert )M^{K}(ds,dx)\right]  \notag \\
& \leq cE\left[ \left( \sum_{t\in J\cap \lbrack 0,T]}((F-F^{\varepsilon
})(t,\Delta Y_{t}^{K}))^{2}\right) ^{1/2}\right]  \notag \\
& \leq c\left( \int_{0}^{T}E\left[ \int \left( \int_{0}^{t}\int \left\vert
(G^{\lambda }-G^{\lambda ,\varepsilon })(x-y)\right\vert
Y_{s}^{K}(dy)ds\right) ^{p}Y_{t}^{K}(dx)\right] dt\right) ^{1/p}  \notag \\
& \rightarrow 0,\;\;\mathrm{as}\;\varepsilon \downarrow 0,\;\forall T>0,
\label{equt:31}
\end{align}%
where the last convergence follows by Proposition \ref{markcorolory} and the
monotone convergence theorem. Now combine (\ref{leonosexqja}), (\ref{equt:5}%
), (\ref{equt:6}), (\ref{condplema}) and (\ref{equt:31}) to get that all the
terms in (\ref{takanaaprox}) converge and the proof of part \textbf{(a)} is
complete. By (\ref{equt:5}), (\ref{equt:6}), (\ref{condplema}) and (\ref%
{equt:31}) we get that all the terms in (\ref{takanaaprox1}) converge and
hence the part \textbf{(b)} of the proposition follows.
\end{proof}

\bigskip

\section{Proof of Proposition~\protect\ref{markcorolory}: estimation of
fractional moments}

\label{moment}


In what follows we will use the following well known equalities. For $p\in
(1,2)$ 
\begin{equation}
z^{p-1}=\eta _{p}\int_{0}^{\infty }\left( 1-e^{-\lambda z}\right) \lambda
^{-p}d\lambda  \label{idepot1}
\end{equation}
and 
\begin{equation}
z^{p}=p\eta _{p}\int_{0}^{\infty }\left( e^{-\lambda z}-1+\lambda z\right)
\lambda ^{-p-1}d\lambda  \label{idepot2}
\end{equation}
where 
\begin{equation*}
\eta _{p}=\frac{p-1}{\Gamma (2-p)}.
\end{equation*}

\begin{proposition}
\label{teomonfrac}Let $1<1+\beta <p<p^{\prime }<2$. Then there exists a
constant $c=c(K,\beta ,p,p^{\prime })$ such that for any non-negative
functions $\varphi ,\psi \in \mathcal{B}_{b}(\mathbb{R}^{d})$, 
\begin{align}
& E\left[ \left. Y_{t}^{K}(\varphi )\left( \int_{0}^{t}Y_{s}^{K}(\psi
)ds\right) ^{p}\right| Y_{0}^{K}=\mu \right]  \notag \\
& \leq c\left\{ \mu \left( S_{t}^{K}\varphi \right) +\mu \left(
S_{t}^{K}\varphi \right) \left( \mu \left( \int_{0}^{t}S_{s}^{K}\psi
ds\right) \right) ^{p}\right.  \notag \\
& +\mu \left( S_{t}^{K}\varphi \right) \mu \left(
\int_{0}^{t}S_{t-s}^{K}\left(\left( \int_{0}^{s}S_{r}^{K}\psi dr\right)
^{p^{\prime }}\right)ds\right)  \notag \\
& +\mu \left( \int_{0}^{t}S_{t-s}^{K}\left( S_{s}^{K}\varphi
\int_{0}^{s}S_{r}^{K}\psi dr\right) ds\right)  \notag \\
& +\mu \left( \int_{0}^{t}S_{t-s}^{K}\left( S_{s}^{K}\varphi
\int_{0}^{s}S_{r}^{K}\psi dr\right) ds\right) \left( \mu \left(
\int_{0}^{t}S_{s}^{K}\psi ds\right) \right) ^{p-1}  \notag \\
& +\mu \left( \int_{0}^{t}S_{t-s}^{K}\left( \int_{0}^{s}S_{s-r}^{K}\left(
S_{r}^{K}\varphi \int_{0}^{r}S_{u}^{K}\psi du\right) dr\left(
\int_{0}^{s}S_{r}^{K}\psi dr\right) ^{p-1}\right) ds\right)  \notag \\
& +\mu \left( \int_{0}^{t}S_{t-s}^{K}\left( S_{s}^{K}\varphi \left(
\int_{0}^{s}S_{r}^{K}\psi dr\right) ^{p}\right) ds\right)  \notag \\
& \left. +\mu \left( \int_{0}^{t}S_{t-s}^{K}\left( S_{s}^{K}\varphi
\int_{0}^{s}S_{s-r}^{K}\left(\left( \int_{0}^{r}S_{u}^{K}\psi du\right)
^{p^{\prime }}\right)dr\right) ds\right) \right\},\;\;\forall t>0.
\label{estfracmom}
\end{align}
\end{proposition}

\begin{proof}
Fix an arbitrary $t>0$. By~(\ref{idepot2}) we obtain 
\begin{eqnarray}
\lefteqn{E_{\mu }\left[ Y_{t}^{K}(\varphi )\left( \int_{0}^{t}Y_{s}^{K}(\psi
)ds\right) ^{p}\right] }  \label{equt:9} \\
&=&p\eta _{p}\int_{0}^{\infty }\lambda ^{-p-1}\left( E\left[
Y_{t}^{K}(\varphi )e^{-\lambda \int_{0}^{t}Y_{s}^{K}(\psi )ds}\right] \right.
\notag \\
&&\left. -E\left[ Y_{t}^{K}(\varphi )\right] +\lambda E\left[
Y_{t}^{K}(\varphi )\int_{0}^{t}Y_{s}^{K}(\psi )ds\right] \right) d\lambda . 
\notag
\end{eqnarray}%
Now we will bound the moments on the right hand side of the above
expression. First of all, by Corollary~\ref{cor:1}, we have 
\begin{eqnarray}
E_{\mu }\left[ Y_{t}^{K}(\varphi )\int_{0}^{t}Y_{s}^{K}(\psi )ds\right]
&=&\mu \left( S_{t}^{K}\varphi \right) \mu \left( \int_{0}^{t}S_{s}^{K}\psi
ds\right)  \label{equt:10} \\
&&+\chi (2)\mu \left( \int_{0}^{t}\int_{0}^{s}S_{r}^{K}(S_{t-r}^{K}\varphi
S_{s-r}^{K}\psi )drds\right) .  \notag
\end{eqnarray}%
Moreover, from Fubini theorem we get the following useful equality 
\begin{equation}
\int_{0}^{t}\int_{0}^{s}S_{r}^{K}(S_{t-r}^{K}\varphi S_{s-r}^{K}\psi
)drds=\int_{0}^{t}S_{t-r}^{K}\left( S_{r}^{K}\varphi
\int_{0}^{r}S_{s}^{K}\psi ds\right) dr.  \label{fubsecmom}
\end{equation}%
Now let us estimate the remaining moment. Use the Laplace transform (\ref%
{lapfuntrun}) and the dominated convergence theorem to obtain 
\begin{align*}
E_{\mu }\left[ Y_{t}^{K}(\varphi )e^{-\lambda \int_{0}^{t}Y_{s}^{K}(\psi )ds}%
\right] & =-\lim_{\varepsilon \downarrow 0}\frac{1}{\varepsilon }E_{\mu }%
\left[ e^{-\lambda \int_{0}^{t}Y_{s}^{K}(\psi )ds-\varepsilon
Y_{t}^{K}(\varphi )}-e^{-\lambda \int_{0}^{t}Y_{s}^{K}(\psi )ds}\right] \\
& =-\lim_{\varepsilon \downarrow 0}\frac{1}{\varepsilon }\left( e^{-\mu
(V_{t}^{K}(\varepsilon \varphi ,\lambda \psi ))}-e^{-\mu
(V_{t}^{K}(0,\lambda \psi ))}\right) .
\end{align*}%
From (\ref{lapfuntrun}), we can easily derive that 
\begin{equation}
V_{s}^{N}(\varphi ,\psi )\geq V_{s}^{N}(0,\psi )\geq 0,\;\;\forall s\geq 0,
\label{equt:11}
\end{equation}%
and hence by the dominated convergence theorem, we get 
\begin{equation*}
V_{t}^{K}(\varepsilon \varphi ,\lambda \psi )\rightarrow V_{t}^{K}(0,\lambda
\psi ),\;\text{as}\;\varepsilon \downarrow 0.
\end{equation*}%
Using the same argument we get 
\begin{equation*}
E\left[ Y_{t}^{K}(\varphi )e^{-\lambda \int_{0}^{t}Y_{s}^{K}(\psi )ds}\right]
=e^{-\mu (V_{t}^{K}(0,\lambda \psi ))}\mu \left( \lim_{\varepsilon
\downarrow 0}\frac{V_{t}^{K}(\varepsilon \varphi ,\lambda \psi
)-V_{t}^{K}(0,\lambda \psi )}{\varepsilon }\right) .
\end{equation*}%
Following the argument in Section 6.3 of \cite{Da2} we can show that $%
U_{t}^{K}(\varphi ,\lambda \psi )$ defined by 
\begin{equation*}
U_{t}^{K}(\varphi ,\lambda \psi )=\lim_{\varepsilon \downarrow 0}\frac{%
V_{t}^{K}(\varepsilon \varphi ,\lambda \psi )-V_{t}^{K}(0,\lambda \psi )}{%
\varepsilon }.
\end{equation*}%
%
%
%
%
%
%
%
%
%
%
satisfies the following equation, 
\begin{equation}
U_{t}^{K}(\varphi ,\lambda \psi )=S_{t}^{K}\varphi
-\int_{0}^{t}S_{t-s}^{K}(U_{s}^{K}(\varphi ,\lambda \psi )(\Phi
^{K})^{\prime }(V_{s}^{K}(0,\lambda \psi )))ds,  \label{newevol}
\end{equation}%
with $\Phi ^{\prime }(x)=\frac{d\Phi (x)}{dx}$ .

Use~(\ref{smfmla}), (\ref{equt:9}), (\ref{equt:10}), (\ref{fubsecmom}) and (%
\ref{newevol}) to get 
\begin{align*}
& E\left[ Y_{t}^{K}(\varphi )\left( \int_{0}^{t}Y_{s}^{K}(\psi )ds\right)
^{p}\right] \\
& =p\eta _{p}\int_{0}^{\infty }\left( e^{-\mu (V_{t}^{K}(0,\lambda \psi
))}\mu (U_{t}^{K}(\varphi ,\lambda \psi ))+\lambda \mu \left(
S_{t}^{K}\varphi \right) \mu \left( \int_{0}^{t}S_{s}^{K}\psi ds\right)
\right. \\
& \left. -\mu \left( S_{t}^{K}\varphi \right) +\lambda \chi (2)\mu \left(
\int_{0}^{t}S_{t-s}^{K}(S_{s}^{K}\varphi \int_{0}^{s}S_{r}^{K}\psi
dr)ds\right) \right) \lambda ^{-p-1}d\lambda \\
& =p\eta _{p}\mu \left( \int_{0}^{\infty }\left( S_{t}^{K}\varphi
e^{-\lambda \mu (\int_{0}^{t}S_{s}^{K}\psi ds)}-S_{t}^{K}\varphi +\lambda
S_{t}^{K}\varphi \mu \left( \int_{0}^{t}S_{s}^{K}\psi ds\right) \right.
\right. \\
& +S_{t}^{K}\varphi e^{-\mu (V_{t}^{K}(0,\lambda \psi ))}-S_{t}^{K}\varphi
e^{-\lambda \mu (\int_{0}^{t}S_{s}^{K}\psi ds)} \\
& +\lambda \chi (2)\int_{0}^{t}S_{t-s}^{K}(S_{s}^{K}\varphi
\int_{0}^{s}S_{r}^{K}\psi dr)ds \\
& \left. \left. -e^{-\mu (V_{t}^{K}(0,\lambda \psi
))}\int_{0}^{t}S_{t-s}^{K}(U_{s}^{K}(\varphi ,\lambda \psi )(\Phi
^{K})^{\prime }(V_{s}^{K}(0,\lambda \psi )))ds\right) \lambda
^{-p-1}d\lambda \right) \\
& =\mu \left( S_{t}^{K}\varphi \right) \left( \mu \left(
\int_{0}^{t}S_{s}^{K}\psi ds\right) \right) ^{p}+I_{1}+I_{2}+I_{3}
\end{align*}%
where 
\begin{eqnarray*}
I_{1} &=&p\eta _{p}\mu \left( S_{t}^{K}\varphi \int_{0}^{\infty }\left(
e^{-\mu (V_{t}^{K}(0,\lambda \psi ))}-e^{-\lambda \mu
(\int_{0}^{t}S_{s}^{K}\psi ds)}\right) \lambda ^{-p-1}d\lambda \right) , \\
I_{2} &=&p\eta _{p}\mu \left( \int_{0}^{\infty }\left( \lambda \chi
(2)\int_{0}^{t}S_{t-s}^{K}(S_{s}^{K}\varphi \int_{0}^{s}S_{r}^{K}\psi
dr)ds\right. \right. \\
&&\left. \left. -\int_{0}^{t}S_{t-s}^{K}(U_{s}^{K}(\varphi ,\lambda \psi
)(\Phi ^{K})^{\prime }(V_{s}^{K}(0,\lambda \psi )))ds\right) \lambda
^{-p-1}d\lambda \right) , \\
I_{3} &=&p\eta _{p}\mu \left( \int_{0}^{\infty }\left( 1-e^{-\mu
(V_{t}^{K}(0,\lambda \psi ))}\right) \right. \\
&&\left. \times \int_{0}^{t}S_{t-s}^{K}(U_{s}^{K}(\varphi ,\lambda \psi
)(\Phi ^{K})^{\prime }(V_{s}^{K}(0,\lambda \psi )))ds\lambda ^{-p-1}d\lambda
\right) .
\end{eqnarray*}

By the elementary inequality $1-e^{-x}\leq x$, for $x\geq 0$, and~(\ref%
{equt:12}) we have 
\begin{equation}
(\Phi ^{K})^{\prime }(x)\leq \eta \frac{K^{1-\beta }}{1-\beta }x= \chi (2)x.
\label{estderivsin}
\end{equation}
Using~(\ref{evoltrunc}) and~(\ref{equt:11}) it is easy to derive that $%
U_{t}^{K}(\varphi ,\lambda \psi )\geq 0$ and 
\begin{eqnarray}  \label{equt:13}
U_{t}^{K}(\varphi ,\lambda \psi )&\leq& S_{t}^{K}\varphi, \\
V_{t}^{K}(0,\lambda \psi )&\leq& \lambda \int_{0}^{t}S_{s}^{K}\psi ds.
\label{estdeuyv}
\end{eqnarray}

The above inequalities and (\ref{newevol}) yield the following bound on $%
I_{3}$: 
\begin{eqnarray}
I_{3} &\leq &c\mu \left( \int_{0}^{\infty
}\int_{0}^{t}S_{t-s}^{K}(S_{s}^{K}\varphi V_{s}^{K}(0,\lambda \psi
))ds\left( 1-e^{-\mu (V_{t}^{K}(0,\lambda \psi ))}\right) \lambda
^{-p-1}d\lambda \right)  \notag \\
&\leq &c\mu \left( \int_{0}^{t}S_{t-s}^{K}(S_{s}^{K}\varphi
\int_{0}^{s}S_{r}^{K}\psi dr)ds\right) \int_{0}^{\infty }\left(
1-e^{-\lambda \mu (\int_{0}^{t}S_{l}^{K}\psi dl)}\right) \lambda
^{-p}d\lambda  \notag \\
&\leq &c\left( \mu \left( \int_{0}^{t}S_{s}^{K}\psi ds\right) \right)
^{p-1}\mu \left( \int_{0}^{t}S_{t-s}^{K}(S_{s}^{K}\varphi
\int_{0}^{s}S_{r}^{K}\psi dr)ds\right) ,\;\;  \label{equt:14}
\end{eqnarray}%
where the last inequality follows by~ (\ref{idepot1}).

Let us take care of $I_{2}$. By~(\ref{newevol}) we get 
\begin{eqnarray}
I_{2} &=&p\eta _{p}\mu \left( \int_{0}^{\infty }\left( \lambda \chi
(2)\int_{0}^{t}S_{t-s}^{K}(S_{s}^{K}\varphi \int_{0}^{s}S_{r}^{K}\psi
dr)ds\right. \right.  \notag  \label{equt:20} \\
&&-\int_{0}^{t}S_{t-s}^{K}(S_{s}^{K}\varphi (\Phi ^{K})^{\prime
}(V_{s}^{K}(0,\lambda \psi )))ds  \notag \\
&&+\int_{0}^{t}S_{t-s}^{K}((\Phi ^{K})^{\prime }(V_{s}^{K}(0,\lambda \psi ))
\notag \\
&&\left. \left. \times \int_{0}^{s}S_{s-r}^{K}(U_{r}^{K}(\varphi ,\lambda
\psi )(\Phi ^{K})^{\prime }(V_{r}^{K}(0,\lambda \psi )))dr)ds\right) \lambda
^{-p-1}d\lambda \right)  \notag \\
&=&J_{1}+J_{2},
\end{eqnarray}%
where 
\begin{eqnarray*}
J_{1} &=&p\eta _{p}\mu \left( \int_{0}^{\infty }\left( \lambda \chi
(2)\int_{0}^{t}S_{t-s}^{K}(S_{s}^{K}\varphi \int_{0}^{s}S_{r}^{K}\psi
dr)ds\right. \right. \\
&&\left. \left. -\int_{0}^{t}S_{t-s}^{K}(S_{s}^{K}\varphi (\Phi
^{K})^{\prime }(V_{s}^{K}(0,\lambda \psi )))ds\right) \lambda
^{-p-1}d\lambda \right) , \\
J_{2} &=&p\eta _{p}\mu \left( \int_{0}^{\infty
}\int_{0}^{t}S_{t-s}^{K}((\Phi ^{K})^{\prime }(V_{s}^{K}(0,\lambda \psi
))\right. \\
&&\left. \times \int_{0}^{s}S_{s-r}^{K}(U_{r}^{K}(\varphi ,\lambda \psi
)(\Phi ^{K})^{\prime }(V_{r}^{K}(0,\lambda \psi )))dr)ds\lambda
^{-p-1}d\lambda \right) .
\end{eqnarray*}

Let us estimate $J_{2}$. First, by~(\ref{defapxintpsim}), (\ref{idepot1})
and (\ref{idepot2}) we obtain 
\begin{eqnarray*}
\int_{0}^{\infty }(\Phi ^{K})^{\prime }(V_{s}^{K}(0,\lambda \psi ))\lambda
^{-p}d\lambda &=&\eta \int_{0}^{\infty }\int_{0}^{K}\left(
1-e^{-wV_{s}^{K}(0,\lambda \psi )}\right) w^{-\beta -1}dw\lambda
^{-p}d\lambda \\
&\leq &\eta \int_{0}^{K}\int_{0}^{\infty }\left( 1-e^{-\lambda
w\int_{0}^{s}S_{r}\psi dr}\right) \lambda ^{-p}d\lambda w^{-\beta -1}dw \\
&=&\eta \eta _{p}^{-1}\int_{0}^{K}\left( w\int_{0}^{s}S_{r}^{K}\psi
dr\right) ^{p-1}w^{-\beta -1}dw \\
&=&c\left( \int_{0}^{s}S_{r}^{K}\psi dr\right) ^{p-1}.
\end{eqnarray*}%
Use this and (\ref{estderivsin}), (\ref{equt:13}), (\ref{estdeuyv}) to get 
\begin{eqnarray}
J_{2} &\leq &c\mu \left( \int_{0}^{\infty }\int_{0}^{t}S_{t-s}^{K}(\right.
(\Phi ^{K})^{\prime }(V_{s}^{K}(0,\lambda \psi ))  \notag \\
&&\left. \times \int_{0}^{s}S_{s-r}^{K}(S_{r}^{K}\varphi V_{r}^{K}(0,\lambda
\psi ))dr)ds\lambda ^{-p-1}d\lambda \right)  \notag \\
&\leq &c\mu \left( \int_{0}^{t}S_{t-s}^{K}\left(
\int_{0}^{s}S_{s-r}^{K}(S_{r}^{K}\varphi \int_{0}^{r}S_{u}^{K}\psi du\right)
dr\right.  \notag \\
&&\left. \times \int_{0}^{\infty }(\Phi ^{K})^{\prime }(V_{s}^{K}(0,\lambda
\psi ))\lambda ^{-p}d\lambda )ds\right)  \notag \\
&\leq &c\mu \left( \int_{0}^{t}S_{t-s}^{K}\left(
\int_{0}^{s}S_{s-r}^{K}\left( S_{r}^{K}\varphi \int_{0}^{r}S_{u}^{K}\psi
du\right) dr\left( \int_{0}^{s}S_{l}^{K}\psi dl\right) ^{p-1}\right)
ds\right) .  \notag \\
&&\mbox{}  \label{equt:15}
\end{eqnarray}

Now let us estimate $J_{1}:$%
\begin{eqnarray*}
J_{1} &=&p\eta _{p}\mu \left( \int_{0}^{\infty }\left( \lambda \chi
(2)\int_{0}^{t}S_{t-s}^{K}\left( S_{s}^{K}\varphi \int_{0}^{s}S_{r}^{K}\psi
dr\right) \right. \right.  \\
&&\left. \left. -\int_{0}^{t}S_{t-s}^{K}\left( S_{s}^{K}\varphi \eta
\int_{0}^{K}(1-e^{-wV_{s}^{K}(0,\lambda \psi )})w^{-\beta -1}dw\right)
ds\right) \lambda ^{-p-1}d\lambda \right)  \\
&=&p\eta _{p}\mu \left( \int_{0}^{t}S_{t-s}^{K}\left( S_{s}^{K}\varphi
\left( \int_{0}^{\infty }\left[ \lambda \chi (2)\int_{0}^{s}S_{r}^{K}\psi
dr\right. \right. \right. \right.  \\
&&\left. \left. \left. +\eta \int_{0}^{K}(e^{-wV_{s}^{K}(0,\lambda \psi
)}-1)w^{-\beta -1}dw\right] \lambda ^{-p-1}d\lambda )\right) ds\right) .
\end{eqnarray*}%
Using the identity 
\begin{equation*}
\chi (2)=\eta \int_{0}^{K}w^{-\beta }dw,
\end{equation*}%
we obtain%
\begin{eqnarray*}
J_{1} &=&p\eta _{p}\mu \left( \int_{0}^{t}S_{t-s}^{K}(S_{s}^{K}\varphi \eta
\int_{0}^{\infty }\int_{0}^{K}(e^{-wV_{s}^{K}(0,\lambda \psi )}-1\right.  \\
&&\left. \left. +\lambda w\int_{0}^{s}S_{r}^{K}\psi dr)w^{-\beta
-1}dw\lambda ^{-p-1}d\lambda \right) ds\right)  \\
&=&p\eta \eta _{p}\mu \left( \int_{0}^{t}S_{t-s}^{K}(S_{s}^{K}\varphi \left[
\int_{0}^{K}\int_{0}^{\infty }(e^{-\lambda w\int_{0}^{s}S_{r}\psi
dr}-1\right. \right.  \\
&&+\lambda w\int_{0}^{s}S_{r}^{K}\psi dr)\lambda ^{-p-1}d\lambda w^{-\beta
-1}dw \\
&&\left. \left. \left. +\int_{0}^{K}\int_{0}^{\infty }\left(
e^{-wV_{s}^{K}(0,\lambda \psi )}-e^{-\lambda w\int_{0}^{s}S_{r}^{K}\psi
dr}\right) \lambda ^{-p-1}d\lambda w^{-\beta -1}dw\right] \right) ds\right) 
\end{eqnarray*}%
\begin{eqnarray}
&=&\mu \left( \int_{0}^{t}S_{t-s}^{K}\left( S_{s}^{K}\varphi \eta
\int_{0}^{K}\left( w\int_{0}^{s}S_{r}^{K}\psi dr\right) ^{p}w^{-\beta
-1}dw\right) \right)   \notag \\
&&+p\eta \eta _{p}\mu \left( \int_{0}^{t}S_{t-s}^{K}(S_{s}^{K}\varphi
\right.   \notag \\
&&\left. \left. \times \int_{0}^{K}\int_{0}^{\infty }\left(
e^{-wV_{s}^{K}(0,\lambda \psi )}-e^{-\lambda w\int_{0}^{s}S_{r}^{K}\psi
dr}\right) \lambda ^{-p-1}d\lambda w^{-\beta -1}dw\right) ds\right)   \notag
\\
&=&c\mu \left( \int_{0}^{t}S_{t-s}^{K}\left( S_{s}^{K}\varphi \left(
\int_{0}^{s}S_{r}^{K}\psi dr\right) ^{p}\right) \right) ds  \notag \\
&&+p\eta \eta _{p}\mu \left( \int_{0}^{t}S_{t-s}^{K}\left( S_{s}^{K}\varphi
Q(s)\right) ds\right) ,
\end{eqnarray}%
where 
\begin{eqnarray*}
Q(s) &=&\int_{0}^{K}\int_{0}^{\infty }(e^{-wV_{s}^{K}(0,\lambda \psi
)}-e^{\lambda w\int_{0}^{s}S_{r}^{K}\psi dr})\lambda ^{-p-1}d\lambda
w^{-\beta -1}dw \\
&\leq &\int_{0}^{K}\int_{0}^{\infty }\left\vert wV_{s}^{K}(0,\lambda \psi
)-\lambda w\int_{0}^{s}S_{r}^{K}\psi dr\right\vert \lambda ^{-p-1}d\lambda
w^{-\beta -1}dw \\
&=&c\left( \int_{0}^{K}+\int_{K}^{\infty }\right) \left\vert
V_{s}^{K}(0,\lambda \psi )-\int_{0}^{s}S_{r}^{K}(\lambda \psi )dr\right\vert
\lambda ^{-p-1}d\lambda  \\
&=&(Q_{1}+Q_{2})(s).
\end{eqnarray*}%
By~(\ref{evoltrunc}) 
\begin{equation*}
Q_{1}(s)=c\int_{0}^{K}\int_{0}^{s}S_{s-r}^{K}\Phi ^{K}(V_{r}^{K}(0,\lambda
\psi ))dr\lambda ^{-p-1}d\lambda .
\end{equation*}%
Due to $1<1+\beta <p<p^{\prime }<2$ we have, from the elementary inequality $%
0\leq e^{-x}-1-x\leq cx^{p^{\prime }},$ for $x\geq 0$, that%
\begin{equation*}
\Phi ^{K}(x)\leq c\eta \frac{K^{p^{\prime }-\beta -1}}{p^{\prime }-\beta -1}%
x^{p^{\prime }}.
\end{equation*}%
Using this we obtain 
\begin{eqnarray}
Q_{1}(s) &\leq &c\int_{0}^{K}\int_{0}^{s}S_{s-r}^{K}\left(
(V_{r}^{K}(0,\lambda \psi ))^{p^{\prime }}\right) dr\lambda ^{-p-1}d\lambda 
\notag  \label{equt:18} \\
&\leq &c\int_{0}^{K}\int_{0}^{s}S_{s-r}^{K}\left( \left(
\int_{0}^{r}S_{u}^{K}(\lambda \psi )du\right) ^{p^{\prime }}\right)
dr\lambda ^{-p-1}d\lambda   \notag \\
&=&c\int_{0}^{s}S_{s-r}^{K}\left( \left( \int_{0}^{r}S_{u}^{K}\psi du\right)
^{p^{\prime }}\right) dr.
\end{eqnarray}%
Apply triangle inequality and (\ref{estdeuyv}) to bound $Q_{2}$: 
\begin{eqnarray}
Q_{2}(s) &=&c\int_{K}^{\infty }\left\vert V_{s}^{K}(0,\lambda \psi
)-\int_{0}^{s}S_{r}^{K}(\lambda \psi )dr\right\vert \lambda ^{-p-1}d\lambda 
\notag  \label{equt:16} \\
&\leq &2c\int_{K}^{\infty }\int_{0}^{s}S_{r}^{K}(\lambda \psi )dr\lambda
^{-p-1}d\lambda   \notag \\
&=&c\int_{0}^{s}S_{r}^{K}\psi dr.
\end{eqnarray}

Finally, let us estimate $I_{1}$. Proceeding as before we have 
\begin{eqnarray}
I_{1} &=&p\eta _{p}\mu \left( S_{t}^{K}\varphi \left(
\int_{0}^{K}+\int_{K}^{\infty }\right) \left( e^{-\mu (V_{t}^{K}(0,\lambda
\psi ))}-e^{-\lambda \mu (\int_{0}^{t}S_{s}^{K}\psi ds)}\right) \lambda
^{-p-1}d\lambda \right)  \notag  \label{equt:19} \\
&\leq &p\eta _{p}\mu \left( S_{t}^{K}\varphi \int_{0}^{K}\left\vert \mu
(V_{t}^{K}(0,\lambda \psi ))-\mu \left( \int_{0}^{t}S_{s}^{K}\psi ds\right)
\right\vert \lambda ^{-p-1}d\lambda \right)  \notag \\
&&+p\eta _{p}\mu \left( S_{t}^{K}\varphi \int_{K}^{\infty }2\lambda
^{-p-1}d\lambda \right)  \notag \\
&\leq &c\mu (S_{t}^{K}\varphi )\int_{0}^{K}\mu \left(
\int_{0}^{t}S_{t-s}^{K}(V_{s}^{K}(0,\lambda \psi ))^{p^{\prime }}ds\right)
\lambda ^{-p-1}d\lambda +c\mu (S_{t}^{K}\varphi )  \notag \\
&=&c\mu (S_{t}^{K}\varphi )\mu \left( \int_{0}^{t}S_{t-s}^{K}\left(
\int_{0}^{s}S_{r}^{K}\psi dr\right) ^{p^{\prime }}ds\right) +c\mu
(S_{t}^{K}\varphi ).
\end{eqnarray}%
Combining (\ref{equt:14})-(\ref{equt:19}) and (\ref{smfmla}) 
we obtain (\ref{estfracmom}).
\end{proof}

\bigskip

Now, the proof of Proposition~\ref{markcorolory} is based on the bounds that
we will get on all the terms on the right hand side of~(\ref{estfracmom}).

\begin{lemma}
\label{lem01} Let $\mu (dx)=h(x)dx$, where $h$ is bounded and integrable.
Then 
\begin{equation}
\sup_{x\in \mathbb{R}^{d}}\mu \left( S_{s}G^{\lambda }(\cdot \cdot
-x)\right) \leq c_{\ref{estpac0}}(h,\lambda)<\infty.  \label{estpac0}
\end{equation}
\end{lemma}

\begin{proof}
Using the explicit expression for $\mu $ we have, 
\begin{eqnarray*}
\mu \left( S_{s}G^{\lambda }(\cdot \cdot -x)\right) &=&\int \int
p_{s}(y-z)G^{\lambda }(z-x)dz\mu (dy) \\
&\leq &\int \int p_{s}(y-z)G^{\lambda }(z-x)dz\left\Vert h\right\Vert
_{\infty }dy \\
&=&\left\Vert h\right\Vert _{\infty }\left\Vert G^{\lambda }\right\Vert
_{1}=\left\Vert h\right\Vert _{\infty }\lambda ^{-1},
\end{eqnarray*}%
recall that $||\cdot ||_{\infty }$ is the supremum norm.
\end{proof}

In the next two lemmas we are going to use the following basic inequalities:
For $d>\alpha $ and $\delta \in (0,1)$, we have (\cite{M}, Lemma 4) 
\begin{equation}
p_{t}(x)\leq ct^{\delta -1}\left\vert x\right\vert ^{\alpha -d-\alpha \delta
}\text{, \ \ \ }t>0,\;x\in \mathbb{R}^{d}\backslash \{0\},
\label{stdensiest}
\end{equation}%
and the Riesz convolution formula 
\begin{equation}
\int_{\mathbb{R}^{d}}\left\vert x-z\right\vert ^{a-d}\left\vert
z-y\right\vert ^{b-d}dz=c\left\vert x-y\right\vert ^{a+b-d},
\label{rieszconv}
\end{equation}%
whenever $a,b>0,a+b<d$ and $x,y\in \mathbb{R}^{d}$.


Also define the indicator function: 
\begin{equation*}
\kappa (x)=\mathbf{1}(\left\vert x\right\vert \leq 1).
\end{equation*}

\begin{lemma}
\label{lem1}Let $\alpha <d$. Then, for any $\delta \in (0,1)$ and $a\in
\lbrack 0,d)$, we have 
\begin{equation}
\int_{\mathbb{R}^{d}}p_{s}(y-z)|z-x|^{-a}\,dz\leq c_{1}+c_{2}s^{\delta
-1}(|y-x|^{\alpha -a-\delta \alpha }\kappa (y-x)+1),\;\;\forall y,x\in 
\mathbb{R}^{d},  \label{desglm1}
\end{equation}%
where $c_{1}\geq 1$, $c_{2}>0$ are constants.
\end{lemma}

\begin{proof}
First let us prove (\ref{desglm1}) for the case $\alpha -\delta \alpha <a$.
Use (\ref{stdensiest}) and (\ref{rieszconv}) to get 
\begin{eqnarray*}
\int_{\mathbb{R}^{d}}p_{s}(y-z)|z-x|^{-a}\,dz &=&\int_{\mathbb{R}%
^{d}}p_{s}(y-x-z)|z|^{-a}\,dz \\
&= &\left( \int_{|z|>1}+\int_{|z|\leq 1}\right) p_{s}(y-x-z)|z|^{-a}\,dz \\
&\leq &1+\int_{|z|\leq 1}s^{\delta -1}|y-x-z|^{\alpha -d-\delta
\alpha}|z|^{-a}\,dz \\
&\leq &1+cs^{\delta -1}|y-x|^{\alpha -a-\delta \alpha } \\
&=&1+cs^{\delta -1}(|y-x|^{\alpha -a-\delta \alpha } \mathbf{1}(\left\vert
y-x\right\vert \leq 1) \\
& &+|y-x|^{\alpha -a-\delta \alpha }\mathbf{1} (\left\vert y-x\right\vert>1))
\\
&\leq &1+cs^{\delta -1}(|y-x|^{\alpha -a-\delta \alpha }\kappa (y-x)+1).
\end{eqnarray*}

Now, suppose $\alpha -\delta \alpha \geq a$. Using a simple coupling
argument, as in Lemma 5.1 of \cite{lea}, we have%
\begin{equation*}
\int_{\mathbb{R}^{d}}p_{s}(y-x-z)|z|^{-a}\,dz\leq \int_{\mathbb{R}%
^{d}}p_{s}(z)|z|^{-a}\,dz.
\end{equation*}%
By the scaling relationship 
\begin{equation}
p_{t}(x)=t^{-d/\alpha }p_{1}(t^{-1/\alpha }x),\text{ \ \ }t>0,\text{ }x\in 
\mathbb{R}^{d},  \label{relaciondeescala}
\end{equation}%
we get%
\begin{eqnarray*}
\int_{\mathbb{R}^{d}}p_{s}(z)|z|^{-a}\,dz &=&s^{-a/\alpha }\int_{\mathbb{R}%
^{d}}p_{1}(s^{-1/\alpha }z)|s^{-1/\alpha }z|^{-a}\,s^{-d/\alpha }dz \\
&=&s^{-a/\alpha }\int_{\mathbb{R}^{d}}p_{1}(z)|z|^{-a}\,dz.
\end{eqnarray*}%
Therefore,%
\begin{eqnarray*}
\int_{\mathbb{R}^{d}}p_{s}(y-z)|z-x|^{-a}\,dz &\leq &cs^{-a/\alpha }(\mathbf{%
1}(\left\vert s\right\vert \leq 1)+\mathbf{1}(\left\vert s\right\vert >1)) \\
&\leq &cs^{\delta -1}+c,
\end{eqnarray*}%
and we are done.
\end{proof}

\begin{lemma}
\label{lem2}For any $\delta \in (0,1)$ there exists $c(t)$ such that for any 
$T>0$, $\sup_{t<T}c(t)<\infty $ and 
\begin{equation}
\int_{0}^{t}S_{s}G^{\lambda }(\cdot -x)(y)ds\leq c(t)(|y-x|^{2\alpha
-d-\delta \alpha }\kappa (y-x)+1),\;\;\forall y,x\in \mathbb{R}^{d}.
\end{equation}
\end{lemma}

\begin{proof}
Let $\alpha <d$. Since $G^{\lambda }(x)\leq c|x|^{\alpha -d}$, take $%
a=d-\alpha $, apply Lemma~\ref{lem1} and make additional integration with
respect to time. If $\alpha \geq d$, then by the unimodality of $p_{1}$,%
\begin{eqnarray*}
\int_{0}^{t}S_{s}G^{\lambda }(\cdot -x)(y)ds &=&\int_{0}^{t}\int_{0}^{\infty
}e^{-\lambda r}p_{r+s}(x-y)\,dr\,ds \\
&\leq &c\int_{0}^{t}\int_{0}^{\infty }e^{-\lambda r}(r+s)^{-d/\alpha
}\,dr\,ds \\
&\leq &c(t),\;\;\forall t\geq 0,
\end{eqnarray*}%
and we are done.
\end{proof}

\begin{lemma}
\label{lem3}Let $1<q<2$, and $d<\alpha (2+1/q)$. Then there exists $c(t)$
such that for any $T>0$, $\sup_{t<T}c(t)<\infty $ and for any $y,x\in 
\mathbb{R}^{d},$ 
\begin{eqnarray}
\int_{0}^{t}S_{t-s}\left( \left( \int_{0}^{s}S_{r}G^{\lambda }(\cdot
-x)\,dr\right) ^{q}\right) (y)ds &\leq &c(t),  \label{equt:25} \\
\int_{0}^{t}S_{s+\varepsilon }\left( \left( \int_{0}^{s}S_{r}G^{\lambda
}(\cdot -x)\,dr\right) ^{q}\right) (y)ds &\leq &c(t),\;\;\forall \varepsilon
\in \lbrack 0,1].
\end{eqnarray}
\end{lemma}

\begin{proof}
Since $d<\alpha (2+1/q)$ it is easy to check that we can fix $\delta \in
(0,1)$ sufficiently small such that, 
\begin{equation}
q(2\alpha -d-\delta \alpha )+\alpha -\delta \alpha >0.  \label{equt:22}
\end{equation}%
%
%
%
%
%
%
%
%
%
%
By Lemma~\ref{lem2}, 
\begin{equation}
\left( \int_{0}^{s}S_{r}G^{\lambda }(\cdot -x)\,dr\right) ^{q}(y)\leq
\sup_{s\leq t}c(s)(|y-x|^{q(2\alpha -d-\delta \alpha )}\kappa (y-x)+1).
\label{equt:23}
\end{equation}%
Now take $a=-q(2\alpha -d-\delta \alpha )$. If $a<0$ then the result follows
trivially, due to the fact that then the right hand side of (\ref{equt:23})
is uniformly bounded for any $y,x\in \mathbb{R}^{d}$.

If $a\geq 0$, we apply again Lemma~\ref{lem1} to conclude that the result
follows if $q(2\alpha -d-\delta \alpha )+\alpha -\delta \alpha >0$. But this
is exactly the condition (\ref{equt:22}) which is satisfied due to the
choice of $\delta $.
\end{proof}

\begin{proof}[Proof of Proposition~\protect\ref{markcorolory}]
From (\ref{relstrunsem}) we see that $S_{t}^{K}\leqslant S_{t}$, hence
Proposition~\ref{teomonfrac} implies 
\begin{align*}
& E_{\mu }\left[ \int_{\mathbb{R}^{d}}Y_{t}^{K}(p_{\varepsilon }(\cdot
-x))\left( \int_{0}^{t}Y_{s}^{K}(G^{\lambda }(\cdot \cdot -x))ds\right)
^{p}dx\right] \\
& \leq c\left\{ \int_{\mathbb{R}^{d}}\mu (S_{t}p_{\varepsilon }(\cdot
-x))dx\right. \\
& +\int_{\mathbb{R}^{d}}\mu (S_{t}p_{\varepsilon }(\cdot -x))\left( \mu
\left( \int_{0}^{t}S_{s}G^{\lambda }(\cdot \cdot -x)ds\right) \right) ^{p}dx
\\
& +\int_{\mathbb{R}^{d}}\mu \left( S_{t}p_{\varepsilon }(\cdot -x)\right)
\mu \left( \int_{0}^{t}S_{t-s}\left( \int_{0}^{s}S_{r}G^{\lambda }(\cdot
\cdot -x)dr\right) ^{p^{\prime }}ds\right) dx \\
& +\int_{\mathbb{R}^{d}}\mu \left( \int_{0}^{t}S_{t-s}\left(
S_{s}p_{\varepsilon }(\cdot -x)\int_{0}^{s}S_{r}G^{\lambda }(\cdot \cdot
-x)dr\right) ds\right) dx \\
& +\int_{\mathbb{R}^{d}}\mu \left( \int_{0}^{t}S_{t-s}\left(
S_{s}p_{\varepsilon }(\cdot -x)\int_{0}^{s}S_{r}G^{\lambda }(\cdot \cdot
-x)dr\right) ds\right) \\
& \times \left( \mu \left( \int_{0}^{t}S_{s}G^{\lambda }(\cdot \cdot
-x)ds\right) \right) ^{p-1}dx
\end{align*}%
\begin{eqnarray*}
&&+\int_{\mathbb{R}^{d}}\mu \left( \int_{0}^{t}S_{t-s}\left(
\int_{0}^{s}S_{s-r}\left( S_{r}p_{\varepsilon }(\cdot
-x)\int_{0}^{r}S_{u}G^{\lambda }(\cdot -x)du\right) dr\right. \right. \\
&&\left. \left. \times \left( \int_{0}^{s}S_{r}G^{\lambda }(\cdot \cdot
-x)dr\right) ^{p-1}\right) ds\right) dx \\
&&+\int_{\mathbb{R}^{d}}\mu \left( \int_{0}^{t}S_{t-s}\left(
S_{s}p_{\varepsilon }(\cdot -x)\left( \int_{0}^{s}S_{r}G^{\lambda }(\cdot
\cdot -x)dr\right) ^{p}\right) ds\right) dx \\
&&\left. +\int_{\mathbb{R}^{d}}\mu \left( \int_{0}^{t}S_{t-s}\left(
S_{s}p_{\varepsilon }(\cdot -x)\int_{0}^{s}S_{s-r}\left(
\int_{0}^{r}S_{u}G^{\lambda }(\cdot \cdot -x)du\right) ^{p^{\prime
}}dr\right) ds\right) dx\right\} \\
&=&c\sum_{i=1}^{8}I_{i}(\varepsilon ).
\end{eqnarray*}%
We will check the boundedness of all the terms $I_{i}(\varepsilon
),\;i=1,\ldots ,8$. First note, that for $d\leq \alpha $ all the terms $%
I_{i}(\varepsilon ),\;i=1,\ldots ,8$ can be bounded very easily, and we
leave it to check to the reader. We will consider the case $\alpha <d$. The
first two terms, $I_{1}(\varepsilon )$ and $I_{2}(\varepsilon )$ are easy to
handle. By the Fubini theorem and Lemma~\ref{lem01} we get 
\begin{equation*}
I_{1}(\varepsilon )+I_{2}(\varepsilon )\leq \mu (1)\left( 1+\left( c_{\ref%
{estpac0}}t\right) ^{p}\right) .
\end{equation*}%
By Lemma~\ref{lem3} we easily get 
\begin{equation*}
I_{3}(\varepsilon )+I_{8}(\varepsilon )\leq \mu (1)c(t).
\end{equation*}%
For $I_{7}(\varepsilon )$ we get the following 
\begin{eqnarray}
I_{7}(\varepsilon ) &=&\int_{0}^{t}\int_{\mathbb{R}^{d}\times \mathbb{R}%
^{d}\times \mathbb{R}^{d}}p_{t-s}(y-z)p_{s+\varepsilon }(z-x)  \notag \\
&&\times \left( \int_{0}^{s}S_{r}G^{\lambda }(z-x)dr\right) ^{p}dsdxdz\mu
(dy)  \notag \\
&=&\int_{0}^{t}\int_{\mathbb{R}^{d}\times \mathbb{R}^{d}}p_{t-s}(y-z)%
\,dzS_{s+\varepsilon }\left( \int_{0}^{s}S_{r}G^{\lambda }(\cdot )dr\right)
^{p}(0)ds\mu (dy)  \notag \\
&=&\mu (1)\int_{0}^{t}S_{s+\varepsilon }\left( \int_{0}^{s}S_{r}G^{\lambda
}(\cdot )dr\right) ^{p}(0)ds  \notag \\
&\leq &\mu (1)c(t),
\end{eqnarray}%
where the last inequality follows by Lemma~\ref{lem3}. It is also easy to
check that 
\begin{eqnarray}
I_{4}(\varepsilon ) &\leq &\int_{\mathbb{R}^{d}}\mu \left(
\int_{0}^{t}S_{t-s}\left( S_{s}p_{\varepsilon }(\cdot -x)\left( \left(
\int_{0}^{s}S_{r}G^{\lambda }(\cdot \cdot -x)dr\right) ^{p}+1\right) \right)
ds\right) dx  \notag \\
&\leq &I_{7}(\varepsilon )+\mu (1)t,  \label{equt:29}
\end{eqnarray}%
and 
\begin{equation*}
I_{5}(\varepsilon )\leq c_{\ref{estpac0}}^{p-1}I_{4}(\varepsilon ).
\end{equation*}%
The last term we have to handle is $I_{6}(\varepsilon ):$ 
\begin{eqnarray*}
I_{6}(\varepsilon ) &\leq &\int_{\mathbb{R}^{3d}}\left(
\int_{x:|y_{1}-x|\leq 1}+\int_{x:|y_{1}-x|>1}\right) \left(
\int_{0}^{t}p_{t-s}(y-y_{1})\right. \\
&&\left. \mbox{}\times \int_{0}^{s}p_{s-r}(y_{1}-z)p_{r+\varepsilon
}(z-x)\int_{0}^{r}S_{u}G^{\lambda }(z-x)dudr\right. \\
&&\left. \times \left( \int_{0}^{s}S_{r}G^{\lambda }(y_{1}-x)dr\right)
^{p-1}ds\right) \,dy_{1}\,dz\,dx\,\mu (dy) \\
&=&I_{6,1}(\varepsilon )+I_{6,2}(\varepsilon ).
\end{eqnarray*}%
Our condition on $d$ implies that we can choose $\delta \in (0,1/3)$
sufficiently small such that 
\begin{equation}
(2\alpha -d)(p+1)-\delta \alpha (p+2)>-d.  \label{equt:33}
\end{equation}%
By (\ref{stdensiest}), Lemma \ref{lem1} and Lemma \ref{lem2} we get%
\begin{eqnarray*}
I_{6,1}(\varepsilon ) &\leq &c(t)\int_{\mathbb{R}^{3d}}\int_{x:|y_{1}-x|\leq
1}\left(
\int_{0}^{t}p_{t-s}(y-y_{1})\int_{0}^{s}p_{s-r}(y_{1}-z)(r+\varepsilon
)^{\delta -1}\right. \\
&&\left. \times \left( |z-x|^{3\alpha -2d-2\delta \alpha }\kappa
(z-x)+|z-x|^{\alpha -d-\delta \alpha }\right) \,dr\right. \\
&&\left. \times \left( |y_{1}-x|^{(2\alpha -d-\delta \alpha )(p-1)}+1\right)
ds\right) \,dy_{1}\,dzdx\,\mu (dy) \\
&=&c(t)\int_{\mathbb{R}^{2d}}\int_{x:|y_{1}-x|\leq 1}\left(
\int_{0}^{t}p_{t-s}(y-y_{1})\int_{0}^{s}(r+\varepsilon )^{\delta -1}\right.
\\
&&\left. \times \int_{\mathbb{R}^{d}}p_{s-r}(y_{1}-z)\left( |z-x|^{3\alpha
-2d-2\delta \alpha }\kappa (z-x)+|z-x|^{\alpha -d-\delta \alpha }\right)
\,dzdr\right. \\
&&\left. \times \left( |y_{1}-x|^{(2\alpha -d-\delta \alpha )(p-1)}+1\right)
ds\right) \,dy_{1}\,dx\,\mu (dy) \\
&\leq &c(t)\int_{\mathbb{R}^{2d}}\int_{x:|y_{1}-x|\leq
1}\int_{0}^{t}p_{t-s}(y-y_{1})\int_{0}^{s}(r+\varepsilon )^{\delta -1} \\
&&\left. \times \left[ 1+(s-r)^{\delta -1}\left( |y_{1}-x|^{4\alpha
-2d-3\delta \alpha }+|y_{1}-x|^{2\alpha -d-2\delta \alpha }+1\right) \right]
dr\right. \\
&&\times \left( |y_{1}-x|^{(2\alpha -d-\delta \alpha )(p-1)}+1\right)
ds\,dy_{1}\,dx\,\mu (dy) \\
&=&c(t)\int_{\mathbb{R}^{2d}}\int_{0}^{t}p_{t-s}(y-y_{1})[\int_{0}^{s}(r+%
\varepsilon )^{\delta -1}dr \\
&&\times \int_{x:|y_{1}-x|\leq 1}(1+|y_{1}-x|^{(2\alpha -d-\delta \alpha
)(p-1)})dx+\int_{0}^{s}(r+\varepsilon )^{\delta -1}(s-r)^{\delta -1}dr \\
&&\times \int_{x:|y_{1}-x|\leq 1}(1+|y_{1}-x|^{(2\alpha -d-\delta \alpha
)(p-1)+4\alpha -2d-3\alpha \delta } \\
&&+|y_{1}-x|^{(2\alpha -d-\delta \alpha )(p-1)+2\alpha -d-2\alpha \delta
}+|y_{1}-x|^{(2\alpha -d-\delta \alpha )(p-1)} \\
&&+|y_{1}-x|^{4\alpha -2d-3\alpha \delta }+|y_{1}-x|^{2\alpha -d-2\alpha
\delta })dx]dy_{1}\,\mu (dy) \\
&\leq &c(t)\mu (1),
\end{eqnarray*}%
where the last inequality follows by (\ref{equt:33}). As for the $%
I_{6,2}(\varepsilon )$, by Lemma~\ref{lem2} we get 
\begin{eqnarray*}
I_{6,2}(\varepsilon ) &\leq &\sup_{x,y_{1}:|y_{1}-x|>1}\left(
\int_{0}^{t}S_{r}G^{\lambda }(y_{1}-x)dr\right) ^{p-1}\int_{\mathbb{R}%
^{3d}}\int_{0}^{t}p_{t-s}(y-y_{1}) \\
&&\mbox{}\times \int_{0}^{s}\int_{0}^{r}\int_{0}^{\infty
}p_{s-r}(y_{1}-z)p_{r+\varepsilon +u+v}(0)e^{-\lambda
v}dv\,du\,dr\,ds\,dy_{1}\,dz\,\mu (dy) \\
&\leq &c(t)\mu (1)\int_{0}^{r}\int_{0}^{s}\int_{0}^{r}\int_{0}^{\infty
}(r+\varepsilon +u+v)^{-d/\alpha }e^{-\lambda v}dv\,du\,dr\,ds
\end{eqnarray*}%
and the last integral is bounded if $d<3\alpha $. By combining all the above
estimates we are done with the first part of the proposition.

Now we are going to prove the second part of the proposition. 
Take 
\begin{equation*}
\varphi (y)=\left( \int_{0}^{t}\int G^{\lambda ,\varepsilon
}(y-x)Y_{s}^{K}(dx)ds\right) ^{p}\text{, \ \ }y\in \mathbb{R}^{d}.
\end{equation*}%
For each $n\in \mathbb{N}$ define the truncations functions, $\varphi
_{n}=\varphi \wedge n$. Then $0\leq \varphi _{n}\uparrow \varphi $, as $%
n\rightarrow 0$. Since $\varphi _{n}$ is bounded and $%
p_{1}(z)dzY_{t}^{K}(dy) $ is a finite measure, we have by the dominated
convergence theorem and the scaling relationship (\ref{relaciondeescala})
the following estimation 
\begin{eqnarray*}
\int \varphi _{n}(x)Y_{t}^{K}(dx) &=&\lim_{\delta \downarrow 0}\int \int
\varphi _{n}(\delta ^{1/\alpha }z+x)p_{1}(z)dzY_{t}^{K}(dx) \\
&\leq &\liminf_{\delta \downarrow 0}\int \int \varphi (\delta ^{1/\alpha
}z+x)p_{1}(z)dzY_{t}^{K}(dx) \\
&=&\liminf_{\delta \downarrow 0}\int \int p_{\delta
}(x-y)Y_{t}^{K}(dx)\varphi (y)dy.
\end{eqnarray*}%
Letting $n\rightarrow \infty $, by the monotone convergence theorem, we have 
\begin{align*}
& \int \left( \int_{0}^{t}Y_{s}^{K}(G^{\lambda ,\varepsilon }(\cdot \cdot
-x))ds\right) ^{p}Y_{t}^{K}(dx) \\
& \leq \liminf_{\delta \downarrow 0}\int \left(
\int_{0}^{t}Y_{s}^{K}(G^{\lambda ,\varepsilon }(\cdot \cdot -x)ds\right)
^{p}Y_{t}^{K}(p_{\delta }(\cdot -x))dx.
\end{align*}%
From the Fatou lemma we get 
\begin{align*}
& E\left[ \int \left( \int_{0}^{t}Y_{s}^{K}(G^{\lambda ,\varepsilon }(\cdot
\cdot -x))ds\right) ^{p}Y_{t}^{K}(dx)\right] \\
& \leq \liminf_{\delta \downarrow 0}E\left[ \int \left(
\int_{0}^{t}Y_{s}^{K}(G^{\lambda ,\varepsilon }(\cdot \cdot -x)ds\right)
^{p}Y_{t}^{K}(p_{\delta }(\cdot -x))dx\right] .
\end{align*}%
Since $G^{\lambda ,\varepsilon }\leq G^{\lambda }$ we have by the already
proven part of Proposition \ref{markcorolory}, 
\begin{align*}
& E\left[ \int \left( \int_{0}^{t}Y_{s}^{K}(G^{\lambda ,\varepsilon }(\cdot
\cdot -x))ds\right) ^{p}Y_{t}^{K}(dx)\right] \\
& \leq \liminf_{\delta \downarrow 0}E\left[ \int \left(
\int_{0}^{t}Y_{s}^{K}(G^{\lambda }(\cdot \cdot -x)ds\right)
^{p}Y_{t}^{K}(p_{\delta }(\cdot -x))dx\right] \\
& \leq c(K,p,d,\alpha ,\beta ).
\end{align*}%
Using once again the monotone convergence theorem, as $\varepsilon
\rightarrow 0,$ we arrive at 
\begin{equation*}
E\left[ \int \left( \int_{0}^{t}Y_{s}^{K}(G^{\lambda }(\cdot \cdot
-x))ds\right) ^{p}Y_{t}^{K}(dx)\right] <c(K,p,d,\alpha ,\beta ),
\end{equation*}%
and we are done.
\end{proof}

\textbf{ACKNOWLEDGEMENTS} \medskip

The second author J. V. would like to express his gratitude for the
opportunity to spend some time as a postdoctoral fellow at the Technion
(Haifa, Israel) where this research was done.

\end{document}